# ROBUST ADAPTIVE IMPORTANCE SAMPLING FOR NORMAL RANDOM VECTORS

BY BENJAMIN JOURDAIN[1] AND JÉRÔME LELONG

*Université Paris-Est and Ecole Nationale Supérieure, Paris*

Adaptive Monte Carlo methods are very efficient techniques designed to tune simulation estimators on-line. In this work, we present an alternative to stochastic approximation to tune the optimal change of measure in the context of importance sampling for normal random vectors. Unlike stochastic approximation, which requires very fine tuning in practice, we propose to use sample average approximation and deterministic optimization techniques to devise a robust and fully automatic variance reduction methodology. The same samples are used in the sample optimization of the importance sampling parameter and in the Monte Carlo computation of the expectation of interest with the optimal measure computed in the previous step. We prove that this highly dependent Monte Carlo estimator is convergent and satisfies a central limit theorem with the optimal limiting variance. Numerical experiments confirm the performance of this estimator: in comparison with the crude Monte Carlo method, the computation time needed to achieve a given precision is divided by a factor between 3 and 15.

**Introduction.** We are interested in the computation of $\mathbb{E}(f(G))$, where $G = (G^1, \dots, G^d)$ is a $d$-dimensional standard normal random vector and $f : \mathbb{R}^d \to \mathbb{R}$ is a measurable function such that $f(G)$ is integrable. This problem is particularly important in mathematical finance, where the calculation of the price and hedging ratios of European options in multidimensional Black–Scholes models amounts to the computation of $\mathbb{E}(f(G))$ for a well-chosen function $f$. The same is true when the underlying assets follow more complex dynamics given by stochastic differential equations, which can be

Received November 2008; revised January 2009.

[1]Supported by the French National Research Agency (ANR) program ADAP'MC and the "Chair Risques Financiers," Fondation du Risque.

discretized using the Euler scheme, for instance. We assume that the random variable $f(G)$ is nonzero and slightly more than square-integrable:

$$\mathbb{P}(f(G)\neq 0)>0, \tag{0.1}$$

$$\forall\theta\in\mathbb{R}^d \qquad \mathbb{E}(f^2(G)e^{-\theta\cdot G})<+\infty. \tag{0.2}$$

By Hölder's inequality, (0.2) holds provided $\mathbb{E}(|f|^{2+\varepsilon}(G))<+\infty$ for some $\varepsilon>0$. For any measurable function $h:\mathbb{R}^d\to\mathbb{R}$ either nonnegative or such that $\mathbb{E}|h(G)|<+\infty$, one has

$$\forall\theta\in\mathbb{R}^d \qquad \mathbb{E}(h(G))=\mathbb{E}(h(G+\theta)e^{-\theta\cdot G-|\theta|^2/2}). \tag{0.3}$$

Applying this equality to $h(x)=f(x)$ and to $h(x)=f^2(x)e^{-\theta\cdot x+|\theta|^2/2}$, one obtains that the expectation and the variance of the random variable $f(G+\theta)e^{-\theta\cdot G-|\theta|^2/2}$ are equal to $\mathbb{E}(f(G))$ and $v^f(\theta)-\mathbb{E}^2(f(G))$, respectively, where

$$v^f(\theta)\stackrel{\text{def}}{=}\mathbb{E}(f^2(G)e^{-\theta\cdot G+|\theta|^2/2}).$$

As a consequence, if $(G_i)_{i\geq 1}$ denotes a sequence of i.i.d. $d$-dimensional standard normal random vectors, for any importance sampling parameter $\theta\in\mathbb{R}^d$,

$$M_n(\theta,f)\stackrel{\text{def}}{=}\frac{1}{n}\sum_{i=1}^{n}f(G_i+\theta)e^{-\theta\cdot G_i-|\theta|^2/2}$$

is an unbiased and convergent estimator of $\mathbb{E}(f(G))$. Since $n\operatorname{Var}(M_n(\theta,f))=v^f(\theta)-\mathbb{E}^2(f(G))$, to improve the accuracy of the estimation for a fixed number $n$ of random samples, one should choose $\theta$ minimizing $v^f(\theta)$. The first section of this paper addresses this minimization problem. First, we check that $v^f$ is a strongly convex function going to infinity at infinity, which ensures the existence of a unique value $\theta^f_\star$ such that $v^f(\theta^f_\star)=\inf_{\theta\in\mathbb{R}^d}v^f(\theta)$. Of course, when $\mathbb{E}(f(G))$ is unknown, in general, so is the function $v^f$. Therefore, direct optimization of this function is not implementable. Using a large deviations argument, Glasserman, Heidelberger and Shahabuddin [8] suggest the use of $\theta$ maximizing $\log|f(\theta)|-\frac{|\theta|^2}{2}$, where, by convention, $\log(0)=-\infty$. However, this choice is not optimal and the numerical search for a local maximum of $\log|f(\theta)|-\frac{|\theta|^2}{2}$ is only possible if the function $f$ possesses some regularity. Under (0.2), the function $v^f$ is infinitely continuously differentiable and such that

$$\nabla_\theta v^f(\theta)=\mathbb{E}((\theta-G)f^2(G)e^{-\theta\cdot G+|\theta|^2/2}). \tag{0.4}$$

At this stage, we can understand the appeal of performing the change of measure (0.3) to transform $\mathbb{E}(f^2(G+\theta)e^{-2\theta\cdot G-|\theta|^2})$ into the above expression of $v^f$: no smoothness assumption on the function $f$ is required in order to differentiate within the expectation.

Arouna [1, 2] takes advantage of the characterization of the optimal parameter $\theta_\star^f$ as the unique solution of the equation $\mathbb{E}((\theta - G)f^2(G)e^{-\theta\cdot G+|\theta|^2/2}) = 0$, in order to approximate it by a Robbins–Monro procedure. The standard Robbins–Monro algorithm explodes, but it can be stabilized using random truncation techniques; see, for instance, [3, 4] or [12]. According to [1], the same random drawings $G_i$ may be used to simultaneously estimate the optimal parameter $\theta_\star^f$ and the expectation of interest $\mathbb{E}(f(G))$. Moreover, both estimators are strongly consistent and the estimator of $\mathbb{E}(f(G))$ is asymptotically normal with an asymptotic variance equal to the optimal one, $v^f(\theta_\star^f) - \mathbb{E}^2(f(G))$. Asymptotic normality of the estimator of $\theta_\star^f$ is discussed in [11].

Tuning the increasing sequence of compact subsets used in randomly truncated procedures is not easy. In [13], Lemaire and Pagès note that using (0.3) in (0.4) leads to $\nabla_\theta v^f(\theta) = e^{|\theta|^2}\mathbb{E}((2\theta - G)f^2(G-\theta))$ and propose to use the characterization of $\theta_\star$ as the unique solution of $\mathbb{E}((2\theta - G)f^2(G-\theta)) = 0$ in order to approximate it by a Robbins–Monro procedure. As soon as the function $f$ satisfies some exponential growth assumptions at infinity, the algorithm they propose is stable without resorting to random truncation techniques. Starting from the present Gaussian framework, Lemaire and Pagès [13] extend this construction of non-exploding Robbins–Monro algorithms to a large class of families of multidimensional probability distributions and even to diffusion process distributions.

In the present paper, we propose and study an alternative approach, one which does not require the delicate tuning of the gain sequence which is still necessary to ensure the stability of Robbins–Monro procedures. When $f(G_i) \neq 0$ for some index $i \in \{1, \dots, n\}$ [by (0.1), a.s. this condition is satisfied for $n$ large enough], the Monte Carlo approximation $v_n^f(\theta) \stackrel{\text{def}}{=} \frac{1}{n}\sum_{i=1}^n f^2(G_i)e^{-\theta\cdot G_i+|\theta|^2/2}$ of the function $v^f$ is also strongly convex and going to infinity at infinity. This ensures the existence of a unique parameter $\theta_n^f$ such that $v_n^f(\theta_n^f) = \inf_{\theta\in\mathbb{R}^d} v_n^f(\theta)$. The function $v_n^f$ is of class $C^\infty$ and its gradient and Hessian matrices

$$\nabla_\theta v_n^f(\theta) = \frac{1}{n}\sum_{i=1}^n (\theta - G_i)f^2(G_i)e^{-\theta\cdot G_i+|\theta|^2/2},$$

$$\nabla_\theta^2 v_n^f(\theta) = \frac{1}{n}\sum_{i=1}^n (I_d + (\theta - G_i)(\theta - G_i)^*)f^2(G_i)e^{-\theta\cdot G_i+|\theta|^2/2}$$

are easily computed if the random samples $(G_i)_{1\le i\le n}$ are stored in the computer memory. Therefore, $\theta_n^f$ can be computed with high precision by a few (four or five) steps of Newton's optimization procedure. In fact, we calculate $\theta_n^f$ as the unique critical point of a modified function $u_n^f$ defined in Section 3.1 and such that the Hessian matrix $\nabla_\theta^2 u_n^f$ is greater than the identity and precised in Section 3.1. In contrast with stochastic approximation procedures, no tuning is necessary for this optimization algorithm. This is the reason for the adjective "robust" in the title of the paper. We propose to estimate $\mathbb{E}(f(G))$ by $M_n(\theta_n^f, f)$. In the context of control variate variance reduction techniques, Kim and Henderson [9] also propose sample average optimization of the control variate parameters as an alternative to stochastic approximation techniques. However, in their algorithm, the expectation of interest is only computed in a second step involving random variables independent of the ones used in the optimization step. In our algorithm, in order to save computation time, the respective approximations $\theta_n^f$ and $M_n(\theta_n^f, f)$ of the optimal parameters and of the expectation of interest are computed using the same set of random samples $(G_i)_{1\le i\le n}$. This makes the mathematical analysis of the properties of $M_n(\theta_n^f, f)$ more complicated: for instance, in general, $M_n(\theta_n^f, f)$ is a biased estimator of $\mathbb{E}(f(G))$. So far, this idea of using the same samples both for the optimization procedure and for the Monte Carlo computation has mainly been investigated in the much simpler context of linear control variates; see [7], Section 4.1 and the references therein. More precisely, for the computation of the expectation $\mathbb{E}(X)$ of the random variable $X$ using $\theta \cdot Y$ as a control variate, where $Y$ is a centered $d$-dimensional random vector, the strong law of large numbers and the central limit theorem enable us to deduce the asymptotic behaviour of the estimator $\frac{1}{n}\sum_{i=1}^n (X_i - \theta_n \cdot Y_i)$, where $(X_i, Y_i)_{i\ge 1}$ is an i.i.d. sample of the law of $(X, Y)$ with empirical mean $(\bar{X}_n, \bar{Y}_n) = \frac{1}{n}\sum_{i=1}^n (X_i, Y_i)$ and $\theta_n = (\sum_{i=1}^n (Y_i - \bar{Y}_n)(Y_i - \bar{Y}_n)^*)^{-1} \sum_{i=1}^n (X_i - \bar{X}_n)(Y_i - \bar{Y}_n)$ minimizes the sample average approximation $\theta \mapsto \frac{1}{n}\sum_{i=1}^n (X_i - \theta \cdot Y_i - \bar{X}_n + \theta \cdot \bar{Y}_n)^2$ of $\mathrm{Var}(X - \theta \cdot Y)$.

The first section of the paper is devoted to the convergence of $\theta_n^f$ to $\theta_\star^f$: almost sure convergence holds and a central limit theorem can be derived under the reinforced integrability condition (1.2). Moreover, $v_n^f(\theta_n^f)$ converges a.s. to $v^f(\theta_\star^f)$. The second section addresses the asymptotic properties, as $n \to \infty$, of our estimator $M_n(\theta_n^f, f)$ of the expectation of interest $\mathbb{E}(f(G))$. We prove that when $f$ is continuous and such that for all $M > 0$, $\mathbb{E}(\sup_{|\theta|\le M} |f(G+\theta)|) < +\infty$, then $M_n(\theta_n^f, f)$ converges a.s. to $\mathbb{E}(f(G))$. In dimension $d = 1$, this continuity assumption may be relaxed: the strong consistency of $M_n(\theta_n^f, f)$ still holds provided $f$ is the sum of a continuous function (as before) and a function of finite variation satisfying some natural growth condition. When $f$ satisfies (1.2) and can be decomposed as the sum of a locally Hölder continuous function with some natural control of the growth of the Hölder continuity constant and of a $C^1$ function satisfying some integrability conditions, then the estimator $M_n(\theta_n^f, f)$ is asymptotically normal with optimal asymptotic variance $v^f(\theta_\star) - \mathbb{E}^2(f(G))$: $\sqrt{n}(M_n(\theta_n^f, f) - \mathbb{E}(f(G))) \xrightarrow{\mathcal{L}} \mathcal{N}_1(0, v^f(\theta_\star) - \mathbb{E}^2(f(G)))$.

Moreover, $\sqrt{n}\frac{M_n(\theta_n^f,f)-\mathbb{E}(f(G))}{\sqrt{v_n^f(\theta_n^f)-M_n^2(\theta_n^f,f)}} \xrightarrow{\mathcal{L}} \mathcal{N}_1(0,1)$, which enables us to construct confidence intervals for $\mathbb{E}(f(G))$. Again, in dimension $d=1$, the conclusion is preserved if one adds a function with finite variation satisfying some natural growth condition to the previous decomposition. In the last section, we illustrate our theoretical results with numerical experiments which confirm the performance of our algorithm.

## 1. Convergence of the importance sampling parameters.

According to our numerical experiments, it may be optimal, in terms of the computation time needed to achieve a given precision for the estimation of $\mathbb{E}(f(G))$, to search for the best importance sampling parameter $\theta$ in a subspace $\{A\vartheta : \vartheta \in \mathbb{R}^{d'}\}$ of $\mathbb{R}^d$, where $A \in \mathbb{R}^{d\times d'}$ is a matrix with rank $d' \le d$. When $f(G)$ corresponds to the payoff of an option written on a $d'$-dimensional Black–Scholes model monitored on a regular time grid, it is sensible to use the same parameter for each coordinate $G^k$ corresponding to a time increment of a given Brownian coordinate. In Section 3.2, the choice of the corresponding matrix $A$ is made precise for the standard Black–Scholes model and the multidimensional Black–Scholes model in the examples of the one-dimensional barrier option and the barrier basket option, respectively. For this choice, according to our numerical experiments, the variances obtained with and without parameter reduction are nearly the same. That is why we introduce

$$v^{f,A}(\vartheta) \stackrel{\text{def}}{=} \mathbb{E}(f^2(G)e^{-A\vartheta\cdot G+|A\vartheta|^2/2}).$$

Since $v^{f,A}(\vartheta) = v^f(A\vartheta)$, the properties of the function $v^{f,A}$ may be deduced from those of $v^f$. The case treated in the Introduction corresponds to the particular choices $d'=d$ and $A=I_d$.

LEMMA 1.1. *Under (0.2), the function $v^f$ is infinitely continuously differentiable with, for all $\alpha=(\alpha^1,\dots,\alpha^d)\in\mathbb{N}^d$ and all $\theta=(\theta^1,\dots,\theta^d)\in\mathbb{R}^d$,*

$$\frac{\partial^{\alpha^1+\cdots+\alpha^d}}{\partial_{\theta^1}^{\alpha^1}\cdots\partial_{\theta^d}^{\alpha^d}}v^f(\theta) = \mathbb{E}\left(\frac{\partial^{\alpha^1+\cdots+\alpha^d}}{\partial_{\theta^1}^{\alpha^1}\cdots\partial_{\theta^d}^{\alpha^d}}[f^2(G)e^{-\theta\cdot G+|\theta|^2/2}]\right).$$

*Under (0.1), the function $v^f$ is strongly convex and hence is such that $\lim_{|\theta|\to+\infty} v^f(\theta) = +\infty$.*

PROOF. The function $\theta \mapsto f^2(G)e^{-\theta\cdot G+|\theta|^2/2}$ is infinitely continuously differentiable with $\frac{\partial}{\partial_{\theta^j}} f^2(G)e^{-\theta\cdot G+|\theta|^2/2} = f^2(G)(\theta^j - G^j)e^{-\theta\cdot G+|\theta|^2/2}$. Since

$$\sup_{|\theta|\le M} |\partial_{\theta^j} f^2(G)e^{-\theta\cdot G+|\theta|^2/2}|$$

(1.1)
$$\leq e^{M^2/2} f^2(G)(M + (e^{G^j} + e^{-G^j})) \prod_{k=1}^{d} (e^{MG^k} + e^{-MG^k}),$$

where the right-hand side is integrable by (0.2), Lebesgue's theorem ensures that $v^f$ is continuously differentiable with $\frac{\partial}{\partial_{\theta j}} v^f(\theta) = \mathbb{E}(f^2(G)(\theta^j - G^j) \times e^{-\theta \cdot G + |\theta|^2/2})$. Higher-order differentiability properties are obtained by similar arguments and, in particular, $\frac{\partial^2}{\partial_{\theta j}\, \partial_{\theta i}} v^f(\theta) = \mathbb{E}((\mathbf{1}_{\{i=j\}} + (\theta^j - G^j)(\theta^i - G^i)) f^2(G) e^{-\theta \cdot G + |\theta|^2/2})$.

Assumption (0.1) ensures the existence of $\varepsilon > 0$ such that $\mathbb{P}(f^2(G) \geq \varepsilon, |G| \leq \frac{1}{\varepsilon}) > 0$. Since $\mathbb{E}(f^2(G)e^{-\theta \cdot G + |\theta|^2/2}) \geq \varepsilon e^{-|\theta|/\varepsilon + |\theta|^2/2} \mathbb{P}(f^2(G) \geq \varepsilon, |G| \leq \frac{1}{\varepsilon})$, one easily deduces that $\lim_{|\theta| \to +\infty} \mathbb{E}(f^2(G)e^{-\theta \cdot G + |\theta|^2/2}) = +\infty$. As the continuous function $\theta \mapsto \mathbb{E}(f^2(G)e^{-\theta \cdot G + |\theta|^2/2})$ does not vanish, the Hessian matrix $\nabla^2_\theta v^f(\theta)$ is uniformly bounded from below by the positive definite matrix $\inf_{\theta \in \mathbb{R}^d} \mathbb{E}(f^2(G)e^{-\theta \cdot G + |\theta|^2/2}) I_d$. This yields the strong convexity of the function $v^f$. □

As a consequence, $v^{f,A}$ is a strongly convex function going to infinity at infinity and there exists a unique $\vartheta^{f,A}_\star \in \mathbb{R}^{d'}$ such that $v^{f,A}(\vartheta^{f,A}_\star) = \inf_{\vartheta \in \mathbb{R}^{d'}} v^{f,A}(\vartheta)$.

Let $(G_i)_{i \geq 1}$ be a sequence of $d$-dimensional independent standard normal random variables. For $n$ large enough, $f(G_i) \neq 0$ for some index $i \in \{1, \dots, n\}$ and the approximation $v^{f,A}_n(\vartheta) = \frac{1}{n} \sum_{i=1}^n f^2(G_i) e^{-A\vartheta \cdot G_i + |A\vartheta|^2/2}$ of $v^{f,A}(\vartheta)$ is also strongly convex and such that $\lim_{|\vartheta| \to +\infty} v^{f,A}_n(\vartheta) = +\infty$. Hence, there exists a unique $\vartheta^{f,A}_n \in \mathbb{R}^{d'}$ such that $v^{f,A}_n(\vartheta^{f,A}_n) = \inf_{\vartheta \in \mathbb{R}^{d'}} v^{f,A}_n(\vartheta)$. The following proposition describes the asymptotic behavior of $\vartheta^{f,A}_n$ as $n \to \infty$. In order to get the central limit theorem, we need the following reinforced integrability condition:

$$\forall \theta \in \mathbb{R}^d \qquad \mathbb{E}(f^4(G)e^{-\theta \cdot G}) < +\infty. \tag{1.2}$$

PROPOSITION 1.2. *Under (0.1) and (0.2), $\vartheta^{f,A}_n$ and $v^{f,A}_n(\vartheta^{f,A}_n)$ converge a.s. to $\vartheta^{f,A}_\star$ and $v^{f,A}(\vartheta^{f,A}_\star)$ as $n \to \infty$. If, moreover, (1.2) holds, then $\sqrt{n}(\vartheta^{f,A}_n - \vartheta^{f,A}_\star) \xrightarrow{\mathcal{L}} \mathcal{N}_{d'}(0, \Gamma)$, where*

$$\Gamma = [\nabla^2 v^{f,A}(\vartheta^{f,A}_\star)]^{-1} \times \operatorname{Cov}(A^*(A\vartheta^{f,A}_\star - G) f^2(G) e^{-A\vartheta^{f,A}_\star \cdot G + |A\vartheta^{f,A}_\star|^2/2}) [\nabla^2 v^{f,A}(\vartheta^{f,A}_\star)]^{-1}$$

*and* $\nabla^2 v^{f,A}(\vartheta) = \mathbb{E}(A^*(I_d + (A\vartheta - G)(A\vartheta - G)^*) A f^2(G) e^{-A\vartheta \cdot G + |A\vartheta|^2/2})$.

REMARK 1.3. The Hessian matrix $\nabla^2 v^{f,A}(\vartheta)$ is positive definite under (0.1) and (0.2). If, moreover, (1.2) holds, using the inequality $|G^k| \le e^{G^k} + e^{-G^k}$ for all $1 \le k \le d$, one obtains that the covariance matrix $\mathrm{Cov}((A\vartheta_\star^{f,A} - G)f^2(G)e^{-A\vartheta_\star^{f,A}\cdot G + |A\vartheta_\star^{f,A}|^2/2})$ exists and $\Gamma$ is well defined.

To obtain some insights into the expression of this asymptotic covariance matrix, note that if $\phi(\vartheta, x) = f^2(x)e^{-A\vartheta\cdot x + |A\vartheta|^2/2}$, then by subtracting $\frac{1}{n} \times \sum_{i=1}^n \nabla_\vartheta \phi(\vartheta_\star^{f,A}, G_i)$ from both sides of the equation $\nabla v_n^{f,A}(\vartheta_n^{f,A}) = \nabla v^{f,A}(\vartheta_\star^{f,A})$ and multiplying by $\sqrt{n}$, one obtains

$$\int_0^1 \frac{1}{n}\sum_{i=1}^n \nabla_\theta^2 \phi(t\vartheta_n^{f,A} + (1-t)\vartheta_\star^{f,A}, G_i)\,dt\sqrt{n}(\vartheta_n^{f,A} - \vartheta_\star^{f,A})$$
$$= \sqrt{n}\left(\mathbb{E}(\nabla_\vartheta \phi(\vartheta_\star^{f,A}, G)) - \frac{1}{n}\sum_{i=1}^n \nabla_\vartheta \phi(\vartheta_\star^{f,A}, G_i)\right).$$

To prove the proposition, we use the following uniform strong law of large numbers, which is a restatement of [14], Lemma A1. This result is also a consequence of the strong law of large numbers in Banach spaces [10], Corollary 7.10, page 189.

LEMMA 1.4. *Let $(X_i)_{i\ge1}$ be a sequence of i.i.d. $\mathbb{R}^m$-valued random vectors and $h:\mathbb{R}^d \times \mathbb{R}^m \to \mathbb{R}$ be a measurable function. Assume that*

- *a.s., $\theta \in \mathbb{R}^d \mapsto h(\theta, X_1)$ is continuous;*
- *$\forall M > 0, \mathbb{E}(\sup_{|\theta|\le M} |h(\theta, X_1)|) < +\infty$.*

*Then, a.s., $\theta \in \mathbb{R}^d \mapsto \frac{1}{n}\sum_{i=1}^n h(\theta, X_i)$ converges locally uniformly to the continuous function $\theta \in \mathbb{R}^d \mapsto \mathbb{E}(h(\theta, X_1))$.*

PROOF OF PROPOSITION 1.2. Since, for $M > 0$,

$$\sup_{|\theta|\le M} f^2(G)e^{-\theta\cdot G + |\theta|^2/2} \le e^{M^2/2} f^2(G)\prod_{k=1}^d (e^{MG^k} + e^{-MG^k}),$$

where the right-hand side is integrable by (0.2), applying Lemma 1.4 with $(X_i)_{i\ge1} = (G_i)_{i\ge1}$ and $h(\theta, x) = f^2(x)e^{-\theta\cdot x + |\theta|^2/2}$ ensures that a.s., $v_n^f$ converges locally uniformly to $v^f$. We restrict ourselves to a subset with probability one of the original probability space on which this convergence holds. Let $\varepsilon > 0$. By the strict convexity and the continuity of $v^{f,A}$,

$$\alpha \stackrel{\text{def}}{=} \inf_{\vartheta:|\vartheta - \vartheta_\star^{f,A}|\ge\varepsilon} v^{f,A}(\vartheta) - v^{f,A}(\vartheta_\star^{f,A}) > 0.$$

The local uniform convergence of $v_n^{f,A}$ to $v^{f,A}$ ensures that

$$\exists n_\alpha \in \mathbb{N}^*, \forall n \ge n_\alpha, \forall \vartheta \text{ s.t. } |\vartheta - \vartheta_\star^{f,A}| \le \varepsilon, \qquad |v_n^{f,A}(\vartheta) - v^{f,A}(\vartheta)| \le \frac{\alpha}{3}.$$

For $n \geq n_\alpha$ and $\vartheta$ such that $|\vartheta - \vartheta_\star^{f,A}| \geq \varepsilon$, we deduce, using the convexity of $v_n^{f,A}$ for the first inequality, that

$$\begin{aligned} &v_n^{f,A}(\vartheta) - v_n^{f,A}(\vartheta_\star^{f,A}) \\ &\quad \geq \frac{|\vartheta - \vartheta_\star^{f,A}|}{\varepsilon}\left[v_n^{f,A}\left(\vartheta_\star^{f,A} + \varepsilon\frac{\vartheta - \vartheta_\star^{f,A}}{|\vartheta - \vartheta_\star^{f,A}|}\right) - v_n^{f,A}(\vartheta_\star^{f,A})\right] \\ &\quad \geq \frac{|\vartheta - \vartheta_\star^{f,A}|}{\varepsilon}\left[v^{f,A}\left(\vartheta_\star^{f,A} + \varepsilon\frac{\vartheta - \vartheta_\star^{f,A}}{|\vartheta - \vartheta_\star^{f,A}|}\right) - v^{f,A}(\vartheta_\star^{f,A}) - \frac{2\alpha}{3}\right] \geq \frac{\alpha}{3}. \end{aligned}$$

Since $v_n^{f,A}(\vartheta_n^{f,A}) \leq v_n^{f,A}(\vartheta_\star^{f,A})$, we conclude that $|\vartheta_n^{f,A} - \vartheta_\star^{f,A}| < \varepsilon$ for $n \geq n_\alpha$. Therefore, $\vartheta_n^{f,A}$ converges a.s. to $\vartheta_\star^{f,A}$. By combining this last result with the local uniform convergence of $v_n^{f,A}$ to the continuous function $v^{f,A}$, we deduce that $v_n^{f,A}(\vartheta_n^{f,A})$ converges a.s. to $v^{f,A}(\vartheta_\star^{f,A})$.

By (1.1) and (0.2), for $M > 0$, $\mathbb{E}(\sup_{|\theta|\leq M} |\nabla_\theta f^2(G)e^{-\theta\cdot G+|\theta|^2/2}|) < +\infty$.

Similarly, $\mathbb{E}(\sup_{|\theta|\leq M} |\nabla_\theta^2 f^2(G)e^{-\theta\cdot G+|\theta|^2/2}|) < +\infty$. The central limit theorem governing the convergence of $\vartheta_n^{f,A}$ to $\vartheta_\star^{f,A}$ follows from [14], Theorem A2. □

## 2. Strong law of large numbers and central limit theorem.

Let

$$\theta_n^{f,A} = A\vartheta_n^{f,A} \quad \text{and} \quad \theta_\star^{f,A} = A\vartheta_\star^{f,A}.$$

The convergence of our estimator $M_n(\theta_n^{f,A}, f)$ of $\mathbb{E}(f(G))$ is ensured by the following theorem, which is a consequence of Propositions 2.6 and 2.13 below. As we do not take advantage of the definition of $\theta_n^{f,A}$ but only use its convergence properties obtained in the previous section, these propositions deal with the asymptotic properties of $M_n(\theta_n^{f,A}, g)$, where $g : \mathbb{R}^d \to \mathbb{R}$ is an arbitrary function and

$$\forall \theta \in \mathbb{R}^d \qquad M_n(\theta, g) \stackrel{\text{def}}{=} \frac{1}{n}\sum_{i=1}^n g(G_i + \theta)e^{-\theta\cdot G_i - |\theta|^2/2}.$$

To make the hypotheses on $f$ precise in the case $d' = 1$ of a one-dimensional importance sampling parameter $\vartheta$, we introduce the following definition.

Definition 2.1. For $\mathcal{A} \in \mathbb{R}^d$, we say that a function $h : \mathbb{R}^d \to \mathbb{R}$:

- is $\mathcal{A}$-*nondecreasing* (resp., $\mathcal{A}$-*nonincreasing*) if

$$\forall x \in \mathbb{R}^d \qquad \vartheta \in \mathbb{R} \mapsto h(x + \mathcal{A}\vartheta) \text{ is nondecreasing (resp., nonincreasing)};$$

- is $\mathcal{A}$-*monotonic* if it is either $\mathcal{A}$-nonincreasing or $\mathcal{A}$-nondecreasing;

- belongs to $\mathcal{V}_{\mathcal{A}}$ if $h$ may be decomposed as the sum of two $\mathcal{A}$-monotonic functions $g_1$ and $g_2$ such that

$$\exists\lambda > 0, \exists\beta \in [0,2), \forall x \in \mathbb{R} \qquad |g_i(x)| \le \lambda e^{|x|^\beta} \quad \text{for } i = 1,2. \tag{2.1}$$

When $d = 1$, $\mathcal{V}_1$ consists of the functions of finite variation which satisfy the growth assumption (2.1). Let us also give an example in a higher dimension.

EXAMPLE. The time-discretization of the one-dimensional Black–Scholes model on the grid $0 = t_0 \le t_1 \le \cdots \le t_d$ is given by $\varphi(G)$, where $\varphi(x) = (e^{(r-\sigma^2/2)t_k + \sigma\sum_{j=1}^k x_j\sqrt{t_j - t_{j-1}}})_{1\le k\le d}$ with $\sigma > 0$. For the choice $A = (\sqrt{t_1}, \sqrt{t_2 - t_1}, \ldots, \sqrt{t_d - t_{d-1}})^*$, which corresponds to the Cameron–Martin formula for the underlying Brownian motion, each coordinate of the function $\varphi$ is $A$-nondecreasing. Therefore, when $g : \mathbb{R}^d \to \mathbb{R}$ is either nondecreasing in each variable or nonincreasing in each variable, the function $g \circ \varphi$ is $A$-monotonic. For $g_1(y) = (y_d - K)^+$ and $g_2(y) = (y_d - K)^+\mathbf{1}_{\{\min_k y_k \ge L\}}$, the functions $g_2 \circ \varphi$ and $(g_1 - g_2) \circ \varphi$, which correspond to the down-and-out and the down-and-in barrier call options also belong to $\mathcal{V}_A$. More generally, all of the barrier call and put option payoffs belong to $\mathcal{V}_A$.

THEOREM 2.2. *Assume (0.1), (0.2) and that $f$ admits a decomposition $f = f_1 + \mathbf{1}_{\{d'=1\}} f_2$, with $f_1$ a continuous function such that $\forall M > 0$, $\mathbb{E}(\sup_{|\theta|\le M} |f_1(G+\theta)|) < +\infty$ and $f_2 \in \mathcal{V}_A$. Then, for any deterministic integer-valued sequence $(\nu_n)_n$ going to $\infty$ with $n$, $M_n(\theta^{f,A}_{\nu_n}, f)$ converges a.s. to $\mathbb{E}(f(G))$.*

Note that for the integrability condition on $f_1$ to hold, it is enough there exist $\beta \in [0,2), \lambda > 0$ such that, for all $x \in \mathbb{R}^d$, $|f_1(x)| \le \lambda e^{|x|^\beta}$.

Under stronger assumptions on $f$, the convergence of $M_n(\theta^{f,A}_n, f)$ to $\mathbb{E}(f(G))$ is governed by a central limit theorem with optimal asymptotic variance $v^{f,A}(\vartheta^{f,A}_\star) - \mathbb{E}^2(f(G))$. For $\alpha \in (0,1]$, let

$$\begin{aligned}\mathcal{H}_\alpha = \{g : \mathbb{R}^d \to \mathbb{R} \text{ s.t. } &\exists\beta \in [0,2), \lambda > 0, \forall x \in \mathbb{R}^d, |g(x)| \le \lambda e^{|x|^\beta} \\ &\forall x, y \in \mathbb{R}^d, |g(x) - g(y)| \le \lambda e^{|x|^\beta \vee |y|^\beta} |x - y|^\alpha\}.\end{aligned}$$

Note that the assumptions of Theorem 2.2 are satisfied for $f \in \mathcal{H}_\alpha$ such that (0.1) holds and that the Hölder condition in the definition of $\mathcal{H}_\alpha$ implies the growth assumption for possibly larger constants $\lambda$ and $\beta$.

THEOREM 2.3. *Assume (0.1), (1.2) and that $f$ admits a decomposition $f = f_1 + f_2 + \mathbf{1}_{\{d'=1\}} f_3$, with $f_1$ a $C^1$ function such that*

$$\forall M > 0 \qquad \mathbb{E}\Big(\sup_{|\theta|\le M} |f_1(G+\theta)| + \sup_{|\theta|\le M} |\nabla f_1(G+\theta)|\Big) < +\infty,$$

*$f_2 \in \mathcal{H}_\alpha$ with $\alpha \in (\frac{\sqrt{d'^2+8d'}-d'}{4}, 1]$ and $f_3 \in \mathcal{V}_A$. Then,*

$$\sqrt{n}(M_n(\theta_n^{f,A}, f) - \mathbb{E}(f(G))) \xrightarrow{\mathcal{L}} \mathcal{N}_1(0, v^{f,A}(\vartheta_\star^{f,A}) - \mathbb{E}^2(f(G))).$$

Note that $\frac{\sqrt{d'^2+8d'}-d'}{4}$ is increasing with $d'$, equals $\frac{1}{2}$ for $d' = 1$ and converges to 1 as $d' \to \infty$. Theorem 2.3 follows from Propositions 2.7 and 2.14 below. With Proposition 1.2, one obtains the following corollary which enables us to construct confidence intervals for $\mathbb{E}(f(G))$ with our algorithm.

Corollary 2.4. *Under the assumptions of Theorem 2.3, if* $\mathrm{Var}(f(G)) > 0$*, then*

$$\sqrt{\frac{n}{v_n^{f,A}(\vartheta_n^{f,A}) - M_n^2(\theta_n^{f,A}, f)}}(M_n(\theta_n^{f,A}, f) - \mathbb{E}(f(G))) \xrightarrow{\mathcal{L}} \mathcal{N}_1(0, 1).$$

Remark 2.5. When $\mathrm{Var}(f(G))$ is positive, the optimal variance $v^{f,A}(\vartheta_\star^{f,A}) - \mathbb{E}^2(f(G))$ is also positive. The estimator $v_n^{f,A}(\vartheta_n^{f,A}) - M_n^2(\theta_n^{f,A}, f)$ converges a.s. to this variance, but may take negative values for $n$ small.

Examples. The hypotheses of Theorems 2.2 and 2.3 are satisfied in the example given after Definition 2.1. Let us give other examples, still inspired by financial applications.

- $f(x) = (K + \sum_{k=1}^d \omega_k e^{\sigma_k (Mx)_k})^+$, where the coefficients $K$, $\omega_k$ and $\sigma_k$ are real numbers and $M \in \mathbb{R}^{d\times d}$: this class of functions belonging to $\mathcal{H}_1$ includes the payoffs of call and put options written on baskets of underlyings in a multidimensional Black–Scholes framework or on a discretely sampled arithmetic average of a single Black–Scholes asset and the payoffs of exchange options on baskets.
- $f(x) = (K + \max_{k=1}^d \omega_k e^{\sigma_k (Mx)_k})^+$, $f(x) = (K + \min_{k=1}^d \omega_k e^{\sigma_k (Mx)_k})^+$: this class of functions belonging to $\mathcal{H}_1$ includes the payoffs of best-of options.
- When $d = 1$, the functions of bounded variation $f(x) = \mathbf{1}_{\{\omega e^{\sigma x} \geq K\}}$ and $f(x) = \mathbf{1}_{\{\omega e^{\sigma x} \leq K\}}$ belong to $\mathcal{V}_1$ and correspond, respectively, to binary call and put options in the Black–Scholes model.
- Let us consider the model

$$dS_t = S_t(\sigma(t, S_t)\, dW_t + r\, dt), \qquad S_0 = s,$$

where $(W_t)_{t\geq 0}$ is a one-dimensional Brownian motion and the local volatility function $\sigma : [0, T] \times \mathbb{R} \to \mathbb{R}$ is bounded and such that $x \mapsto x\sigma(t, x)$ is Lipschitz continuous uniformly for $t \in [0, T]$. When discretizing this SDE by the Euler scheme with $d$ steps of length $h = T/d$ on $[0, T]$, one approximates $S_T$ by $\varphi(G)$, where $\varphi(x) = \phi_d(x_d, \phi_{d-1}(x_{d-1}, \ldots, \phi_1(x_1, s)))$ with $\phi_k(u,v) = v(1+\sigma((k-1)h,v)\sqrt{h}u + rh)$, and $G = \frac{1}{\sqrt{h}}(W_h, W_{2h} - W_h, \ldots, W_{dh} - W_{(d-1)h})$.

There exists $C > 0$ such that, for all $k \in \{1,\ldots,d\}$,

$$\forall u,v,u',v' \in \mathbb{R} \qquad |\phi_k(u,v)| \le C|v|(1+|u|),$$
$$|\phi_k(u,v) - \phi_k(u',v')| \le C((1+(|u| \vee |u'|))|v-v'| + (|v| \vee |v'|)|u-u'|).$$

One deduces, by induction, that for $x,y \in \mathbb{R}^d$, $|\varphi(x)| \le C^d|s|\prod_{k=1}^d(1+|x_k|) \le C^d|s|e^{\sqrt{d}|x|}$ and

$$|\varphi(x)-\varphi(y)| \le C^d|s|\sum_{k=1}^d |x_k - y_k| \prod_{\substack{j=1\\ j\ne k}}^d (1+(|x_j| \vee |y_j|))$$
$$\le C^d|s|\sqrt{d}e^{d(|x|\vee|y|)}|x-y|.$$

Hence, the functions $f(x) = (\varphi(x) - K)^+$ and $f(x) = (K - \varphi(x))^+$ corresponding to the call and put payoffs in the discretized model belong to $\mathcal{H}_1$.

We are now going to study the convergence properties of $M_n(\theta_n^{f,A}, g)$ in the multidimensional framework $d' \ge 1$ before obtaining stronger results in the case $d' = 1$ of a one-dimensional importance sampling parameter.

2.1. *The general case.*

Proposition 2.6. *Let $(\theta_n)_{n\ge1}$ be a sequence of $d$-dimensional random vectors converging almost surely to some random vector $\theta_\infty$ and $g : \mathbb{R}^d \to \mathbb{R}$ be a continuous function such that $\forall M > 0, \mathbb{E}(\sup_{|\theta|\le M} |g(G+\theta)|) < \infty$. Then, $M_n(\theta_n, g)$ converges a.s. to $\mathbb{E}(g(G))$.*

Proof. We apply Lemma 1.4 with $(X_i)_{i\ge1} = (G_i)_{i\ge1}$ and $h(\theta,x) = g(x+\theta)e^{-\theta\cdot x - |\theta|^2/2}$. The continuity assumption follows from the continuity of $g$. Concerning the integrability condition, we deduce from (0.3) and the inequality

$$\sup_{|\theta|\le M} (|g(G+\theta)|e^{-\theta\cdot G - |\theta|^2/2}) \le \sup_{|\theta|\le M} |g(G+\theta)| \prod_{k=1}^d (e^{MG^k} + e^{-MG^k})$$

that

$$\mathbb{E}\Big(\sup_{|\theta|\le M} (|g(G+\theta)|e^{-\theta\cdot G-|\theta|^2/2})\Big) \le e^{dM^2/2} \sum_{\mu\in\{-M,M\}^d} \mathbb{E}\Big(\sup_{|\theta|\le M} |g(G+\theta+\mu)|\Big)$$
$$\le 2^d e^{dM^2/2}\mathbb{E}\Big(\sup_{|\theta|\le(1+\sqrt{d})M} |g(G+\theta)|\Big).$$

Therefore, a.s., $\theta \mapsto M_n(g,\theta)$ converges locally uniformly to the constant function $\theta \mapsto \mathbb{E}(h(\theta,G)) = \mathbb{E}(g(G))$. We easily conclude with the a.s. convergence of $\theta_n$ to $\theta_\infty$. $\square$

PROPOSITION 2.7. *Assume that* $g:\mathbb{R}^d \to \mathbb{R}$ *is such that* $\mathbb{E}(g^2(G+\theta_\star^{f,A}) \times e^{-2\theta_\star^{f,A}\cdot G}) < +\infty$ *and admits a decomposition* $g = g_1 + g_2$, *with* $g_1$ *of class* $C^1$*and satisfying*

$$\forall M > 0 \qquad \mathbb{E}\Big(\sup_{|\theta|\le M} |g_1(\theta+G)| + \sup_{|\theta|\le M} |\nabla g_1(\theta+G)|\Big) < \infty, \tag{2.2}$$

*and* $g_2 \in \mathcal{H}_\alpha$ *for* $\alpha \in (\frac{\sqrt{d'^2+8d'}-d'}{4}, 1]$. *Then, under (0.1) and (1.2),*

$$\sqrt{n}(M_n(\theta_n^{f,A},g) - \mathbb{E}(g(G))) \xrightarrow{\mathcal{L}} \mathcal{N}_1(0, \mathrm{Var}(g(G+\theta_\star^{f,A})e^{-\theta_\star^{f,A}\cdot G - |\theta_\star^{f,A}|^2/2})).$$

By the central limit theorem, $\sqrt{n}(M_n(\theta_\star^{f,A},g) - \mathbb{E}(g(G))) \xrightarrow{\mathcal{L}} \mathcal{N}_1(0, \mathrm{Var}(g(G+\theta_\star^{f,A})e^{-\theta_\star^{f,A}\cdot G - |\theta_\star^{f,A}|^2/2}))$. As a consequence, it is enough to check that for $i \in \{1,2\}$, $\sqrt{n}(M_n(\theta_n^{f,A},g_i) - M_n(\theta_\star^{f,A},g_i)) \xrightarrow{\mathrm{Pr}} 0$. The next lemma deals with the case $i = 1$.

LEMMA 2.8. *Let* $g:\mathbb{R}^d \longrightarrow \mathbb{R}$ *be a* $C^1$ *function satisfying (2.2). Then, under (0.1) and (1.2),* $\sqrt{n}(M_n(\theta_n^{f,A},g) - M_n(\theta_\star^{f,A},g)) \xrightarrow{\mathrm{Pr}} 0$.

Since, for $\varepsilon > 0$,

$$\begin{aligned} &\mathbb{P}(\sqrt{n}|M_n(\theta_n^{f,A},g_2) - M_n(\theta_\star^{f,A},g_2)| \ge \varepsilon) \\ &\quad \le \mathbb{P}(n^\delta|\vartheta_n^{f,A} - \vartheta_\star^{f,A}| \ge 1) \\ &\qquad + \mathbb{P}\Big(\sup_{|\vartheta-\vartheta_\star^{f,A}|\le 1/n^\delta} \sqrt{n}|M_n(A\vartheta,g_2) - M_n(A\vartheta_\star^{f,A},g_2)| \ge \varepsilon\Big), \end{aligned}$$

choosing $\delta \in (d'/2\alpha(d'+2\alpha), 1/2)$, which is possible since $\alpha > \frac{\sqrt{d'^2+8d'}-d'}{4}$, the case $i = 2$ follows from the central limit theorem governing the convergence of $\vartheta_n^{f,A}$ to $\vartheta_\star^{f,A}$ (see Proposition 1.2) combined with the following result.

PROPOSITION 2.9. *Letting* $\mathcal{A} \in \mathbb{R}^{d\times d'}$ *and* $g \in \mathcal{H}_\alpha$ *for* $\alpha \in (0,1]$,

$$\forall \delta > \frac{d'}{2\alpha(d'+2\alpha)}, \forall \vartheta_0 \in \mathbb{R}^{d'}$$

$$\sup_{|\vartheta-\vartheta_0|\le 1/n^\delta} \sqrt{n}|M_n(\mathcal{A}\vartheta,g) - M_n(\mathcal{A}\vartheta_0,g)| \xrightarrow{\mathrm{Pr}} 0.$$

REMARK 2.10. By the argument given in the case $i=2$, if $g \in \mathcal{H}_\alpha$ for $\alpha > \frac{\sqrt{d'^2+8d'}-d'}{4}$, then

$$\sqrt{n}(M_n(\theta_{\nu_n}^{f,A}, g) - \mathbb{E}(g(G))) \xrightarrow{\mathcal{L}} \mathcal{N}_1(0, \operatorname{Var}(g(G+\theta_\star^{f,A})e^{-\theta_\star^{f,A}\cdot G - |\theta_\star^{f,A}|^2/2}))$$

for any deterministic integer-valued sequence $(\nu_n)_n$ such that $\exists \lambda > 0, \exists \gamma > \frac{d'}{\alpha(d'+2\alpha)}, \forall n \in \mathbb{N}^*, \nu_n \geq \lambda n^\gamma$.

PROOF OF LEMMA 2.8. The function $\theta \mapsto M_n(\cdot, g)$ is of class $C^1$ and it is easy to check that $\nabla_\theta M_n(\theta, g) = M_n(\theta, \bar{g})$ with $\bar{g}(x) = \nabla g(x) - g(x)x$. The mean value theorem gives $\sqrt{n}(M_n(\theta_n^{f,A}, g) - M_n(\theta_\star^{f,A}, g)) = A\sqrt{n}(\vartheta_n^{f,A} - \vartheta_\star^{f,A}) \cdot M_n(\bar{\theta}_n^{f,A}, \bar{g})$, with $\bar{\theta}_n^{f,A} \in (\theta_n^{f,A}, \theta_\star^{f,A})$. Since, by Proposition 1.2, $\sqrt{n}(\vartheta_n^{f,A} - \vartheta_\star^{f,A})$ converges in law to a normal random variable, it is enough to prove that $M_n(\bar{\theta}_n^{f,A}, \bar{g}) \xrightarrow{\Pr} 0$. The a.s. convergence of $\vartheta_n^{f,A}$ to $\vartheta_\star^{f,A}$ implies the a.s. convergence of $\bar{\theta}_n^{f,A}$ to $\theta_\star^{f,A}$. Since

$$\sup_{|\theta| \leq M} |(G+\theta)g(G+\theta)| \leq \left(\sum_{k=1}^d (e^{G^k} + e^{-G^k}) + M\right) \sup_{|\theta| \leq M} |g(G+\theta)|,$$

(2.2) combined with the reasoning used at the beginning of the proof of Proposition 2.6 yields

$$\forall M > 0 \qquad \mathbb{E}\Big(\sup_{|\theta| \leq M} |(G+\theta)g(G+\theta)| + \sup_{|\theta| \leq M} |\nabla g(G+\theta)|\Big) < +\infty. \tag{2.3}$$

Proposition 2.6 then implies that $M_n(\bar{\theta}_n^{f,A}, \bar{g}) \xrightarrow{\Pr} \mathbb{E}(\bar{g}(G))$. By (2.3) and the reasoning used at the beginning of the proof of Proposition 2.6,

$$\forall M > 0 \qquad \mathbb{E}\Big(\sup_{|\theta| \leq M} |\bar{g}(G+\theta)| e^{-\theta\cdot G - |\theta|^2/2}\Big) < +\infty.$$

Hence, Lebesgue's theorem implies that $\nabla_\theta \mathbb{E}(g(G+\theta)e^{-\theta\cdot G-|\theta|^2/2}) = \mathbb{E}(\bar{g}(G+\theta)e^{-\theta\cdot G - |\theta|^2/2})$. Since the left-hand side is equal to 0, one deduces for $\theta = 0$ that $\mathbb{E}(\bar{g}(G)) = 0$. $\square$

REMARK 2.11. Let $g : \mathbb{R}^d \to \mathbb{R}$ be a $C^2$ function, $\bar{g}(x) = \nabla g(x) - g(x)x$ and $\bar{\bar{g}}(x) \stackrel{\text{def}}{=} \nabla^2 g(x) - g(x)I_d - x\nabla^* g(x) - \nabla g(x)x^* + g(x)xx^*$. Assume that

$$\mathbb{E}((|g(\theta_\star^{f,A}+G)|^2 + |\bar{g}(\theta_\star^{f,A}+G)|^2)e^{-2\theta_\star^{f,A}\cdot G}) < +\infty, \tag{2.4}$$

$$\begin{aligned} \forall M > 0 \qquad & \mathbb{E}\Big(\sup_{|\theta| \leq M} |g(\theta+G)| \\ & + \sup_{|\theta| \leq M} |\nabla g(\theta+G)| + \sup_{|\theta| \leq M} |\nabla^2 g(\theta+G)|\Big) < \infty. \end{aligned} \tag{2.5}$$

Let $(\nu_n)_n$ be a deterministic integer-valued sequence such that $\exists\lambda>0$, $\forall n\in\mathbb{N}^*$, $\nu_n\ge\lambda\sqrt{n}$. Then, using the decomposition

$$\begin{aligned}\sqrt{n}&(M_n(\theta^{f,A}_{\nu_n},g)-M_n(\theta^{f,A}_\star,g))\\&=\frac{1}{\sqrt{\nu_n}}\sqrt{n}M_n(\theta^{f,A}_\star,\bar g)\cdot\sqrt{\nu_n}(\theta^{f,A}_{\nu_n}-\theta^{f,A}_\star)\\&\quad+\frac{\sqrt{n}}{\nu_n}\sqrt{\nu_n}(\theta^{f,A}_{\nu_n}-\theta^{f,A}_\star)^*\left(\int_0^1(1-t)M_n(t\theta^{f,A}_{\nu_n}+(1-t)\theta^{f,A}_\star,\bar{\bar g})\,dt\right)\\&\quad\times\sqrt{\nu_n}(\theta^{f,A}_{\nu_n}-\theta^{f,A}_\star),\end{aligned}$$

one obtains that under (0.1) and (1.2), the left-hand side converges in probability to 0. As a consequence, $\sqrt{n}(M_n(\theta^{f,A}_{\nu_n},g)-\mathbb{E}(g(G)))\xrightarrow{\mathcal{L}}\mathcal{N}_1(0,\mathrm{Var}(g(G+\theta^{f,A}_\star)e^{-\theta^{f,A}_\star\cdot G-|\theta^{f,A}_\star|^2/2}))$. More generally, if $g$ is of class $C^k$ and satisfies moment assumptions like (2.4) and (2.5) involving its derivatives up to order $k-1$ and $k$, respectively, this result is preserved if $\exists\lambda>0,\forall n\in\mathbb{N}^*,\nu_n\ge\lambda n^{1/k}$.

In order to prove Proposition 2.9, we need the following lemma.

LEMMA 2.12. *If $g\in\mathcal{H}_\alpha$ for $\alpha\in(0,1]$, then*

$$\forall M>0,\exists C>0,\forall\theta,\theta'\in\bar B(0,M),\forall n\in\mathbb{N}^*$$
$$\mathbb{E}((M_n(\theta,g)-M_n(\theta',g))^2)\le\frac{C|\theta-\theta'|^{2\alpha}}{n}.$$

PROOF. Let $M>0$. Since, by (0.3),

$$\mathbb{E}(g(G+\theta)e^{-\theta G-\theta^2/2}-g(G+\theta')e^{-\theta' G-\theta'^2/2})=0,$$

it is enough to check that

$$\exists C>0,\forall\theta,\theta'\in\bar B(0,M)$$
$$\mathbb{E}((g(G+\theta)e^{-\theta\cdot G-|\theta|^2/2}-g(G+\theta')e^{-\theta'\cdot G-|\theta'|^2/2})^2)\le C|\theta-\theta'|^{2\alpha}.$$

One has

$$\begin{aligned}&\mathbb{E}((g(G+\theta)e^{-\theta\cdot G-|\theta|^2/2}-g(G+\theta')e^{-\theta'\cdot G-|\theta'|^2/2})^2)\\&\quad\le2\mathbb{E}((g(G+\theta)-g(G+\theta'))^2e^{-2\theta\cdot G-|\theta|^2})\\&\qquad+2\mathbb{E}(g^2(G+\theta')(e^{-\theta\cdot G-|\theta|^2/2}-e^{-\theta'\cdot G-|\theta'|^2/2})^2).\end{aligned}$$

Let $\lambda>0$ and $\beta\in[0,2)$ be such that

$$\forall x\in\mathbb{R}^d\qquad|g(x)|\le\lambda e^{|x|^\beta},\tag{2.6}$$

$$\forall x,y\in\mathbb{R}^d\qquad|g(x)-g(y)|\le\lambda e^{|x|^\beta\vee|y|^\beta}|x-y|^\alpha.\tag{2.7}$$

One has

$$\text{for } c = 2^{(\beta-1)^+}, \forall a, b \geq 0 \qquad (a+b)^\beta \leq c(a^\beta + b^\beta). \tag{2.8}$$

Since, for $\theta \in \bar{B}(0, M)$, $|\nabla_\theta e^{-\theta\cdot G - \theta^2/2}| = |(G+\theta)e^{-\theta\cdot G-\theta^2/2}| \leq (|G|+M)e^{M|G|}$, one deduces that, for $\theta, \theta' \in \bar{B}(0, M)$,

$$\begin{aligned}
&\mathbb{E}((g(G+\theta)e^{-\theta\cdot G - |\theta|^2/2} - g(G+\theta')e^{-\theta'\cdot G - |\theta'|^2/2})^2) \\
&\quad \leq 2\lambda^2 e^{2cM^\beta}\mathbb{E}((|\theta-\theta'|^{2\alpha} + |\theta-\theta'|^2(|G|+M)^2)e^{2M|G|+2c|G|^\beta}) \\
&\quad \leq C|\theta-\theta'|^{2\alpha}.
\end{aligned}$$

□

Proof of Proposition 2.9. Let $\varepsilon > 0$.

$$\begin{aligned}
&\mathbb{P}\Big(\sup_{|\vartheta-\vartheta_0|\leq 1/n^\delta} \sqrt{n}|M_n(\mathcal{A}\vartheta, g) - M_n(\mathcal{A}\vartheta_0, g)| > \varepsilon\Big) \\
&\qquad \leq n\mathbb{P}(|G| > \sqrt{2d\log n}) \\
&\qquad\quad + \mathbb{P}\Big(\sup_{|\vartheta-\vartheta_0|\leq 1/n^\delta} \sqrt{n}|M_n(\mathcal{A}\vartheta, g) - M_n(\mathcal{A}\vartheta_0, g)| > \varepsilon, \\
&\qquad\qquad\qquad \max_{1\leq i\leq n} |G_i| \leq \sqrt{2d\log n}\Big).
\end{aligned}$$

Since

$$\begin{aligned}
&\mathbb{P}(|G| > \sqrt{2d\log n}) \\
&\quad \leq \sum_{k=1}^d \mathbb{P}(|G^k| > \sqrt{2\log n}) \\
&\quad = 2d\mathbb{P}(G^1 > \sqrt{2\log n}) \leq \frac{2d}{\sqrt{2\log n}} e^{-(\sqrt{2\log n})^2/2},
\end{aligned} \tag{2.9}$$

the second term of the right-hand side tends to 0 as $n$ goes to infinity. Now, let us focus on the first term.

Let $M = |\vartheta_0| + 1$ and $\tilde{M} = |\mathcal{A}|M$. For $\vartheta', \vartheta \in \bar{B}(0, M)$, using (2.7), (2.8) and (2.6) for the second inequality, one obtains

$$\begin{aligned}
&|M_n(\mathcal{A}\vartheta', g) - M_n(\mathcal{A}\vartheta, g)| \\
&\quad \leq \frac{1}{n}\sum_{i=1}^n |g(G_i + \mathcal{A}\vartheta') - g(G_i + \mathcal{A}\vartheta)|e^{-\mathcal{A}\vartheta'\cdot G_i - |\mathcal{A}\vartheta'|^2/2} \\
&\qquad + \frac{1}{n}\sum_{i=1}^n |g(G_i + \mathcal{A}\vartheta)||e^{-\mathcal{A}\vartheta'\cdot G_i - |\mathcal{A}\vartheta'|^2/2} - e^{-\mathcal{A}\vartheta\cdot G_i - |\mathcal{A}\vartheta|^2/2}| \\
&\quad \leq \frac{\lambda|\mathcal{A}|^\alpha|\vartheta'-\vartheta|^\alpha}{n}\sum_{i=1}^n e^{c(|G_i|^\beta + \tilde{M}^\beta) + \tilde{M}|G_i|}(1 + (2\tilde{M})^{1-\alpha}(\tilde{M} + |G_i|)).
\end{aligned}$$

Hence, when $\max_{1\le i\le n}|G_i| \le \sqrt{2d\log n}$, there exists a constant $\gamma$ not depending on $n$ such that if $\nu \stackrel{\text{def}}{=} \frac{\beta\vee 1}{2}$,

$$\begin{aligned}&\forall \vartheta', \vartheta \in \bar{B}(0,M)\\ &\qquad |M_n(\mathcal{A}\vartheta', g) - M_n(\mathcal{A}\vartheta, g)| \le \gamma|\vartheta-\vartheta'|^{\alpha} e^{\gamma(\log n)^{\nu}}.\end{aligned} \tag{2.10}$$

We can cover $\bar{B}(\vartheta_0, \frac{1}{n^\delta})$ with $K = \mathcal{C}\lceil(\gamma^{1/\alpha} n^{1/(2\alpha)-\delta} e^{(\gamma/\alpha)(\log n)^{\nu}}/\varepsilon^{1/\alpha})^{d'}\rceil$ balls of radius $(\frac{\varepsilon}{2\gamma e^{\gamma(\log n)^{\nu}}\sqrt{n}})^{1/\alpha}$, where $\mathcal{C}$ is a geometrical constant not depending on $n$. For $k \in \{1,\dots,K\}$, let $B_k$ denote the $k$th ball and $\vartheta_k$ its center. By (2.10), when $\max_{1\le i\le n}|G_i| \le \sqrt{2d\log n}$,

$$\forall k \in \{1,\dots,K\} \qquad \sup_{\vartheta\in B_k} |M_n(\mathcal{A}\vartheta, g) - M_n(\mathcal{A}\vartheta_k, g)| \le \frac{\varepsilon}{2\sqrt{n}}.$$

Using Lemma 2.12 for the fourth inequality, one deduces that

$$\begin{aligned}&\mathbb{P}\Big(\sup_{|\vartheta-\vartheta_0|\le 1/n^{\delta}} \sqrt{n}|M_n(\mathcal{A}\vartheta, g) - M_n(\mathcal{A}\vartheta_0, g)| > \varepsilon, \max_{1\le i\le n}|G_i| \le \sqrt{2d\log n}\Big)\\
&\quad\le \mathbb{P}\Big(\exists k \le K : |M_n(\mathcal{A}\vartheta_k, g) - M_n(\mathcal{A}\vartheta_0, g)|\\
&\qquad\qquad > \frac{\varepsilon}{\sqrt{n}} - \sup_{\vartheta\in B_k}|M_n(\mathcal{A}\vartheta, g) - M_n(\mathcal{A}\vartheta_k, g)|,\\
&\qquad\qquad\qquad \max_{1\le i\le n}|G_i| \le \sqrt{2d\log n}\Big)\\
&\quad\le \mathbb{P}\Big(\max_{k\le K}|M_n(\mathcal{A}\vartheta_k, g) - M_n(\mathcal{A}\vartheta_0, g)| > \frac{\varepsilon}{2\sqrt{n}}\Big)\\
&\quad\le \sum_{k\le K}\frac{4n}{\varepsilon^2}\mathbb{E}((M_n(\mathcal{A}\vartheta_k, g) - M_n(\mathcal{A}\vartheta_0, g))^2)\end{aligned}$$

$$\le \sum_{k\le K}\frac{4n}{\varepsilon^2}\frac{C|\vartheta-\vartheta_k|^{2\alpha}}{n} \le Cn^{d'/(2\alpha)-(d'+2\alpha)\delta} e^{d'\gamma/\alpha(\log n)^{\nu}}. \tag{2.11}$$

Since $\beta < 2$ and $\delta > \frac{d'}{2\alpha(d'+2\alpha)}$, $\nu < 1$ and $\frac{d'}{2\alpha} - (d'+2\alpha)\delta < 0$. Therefore, the upper bound in equation (2.11) converges to 0 as $n$ increases to infinity. □

2.2. *The case of a one-dimensional importance sampling parameter.* In the present section, dedicated to the case $d' = 1$ of a one-dimensional importance sampling parameter, we obtain convergence results under weaker assumptions on the function $g$.

PROPOSITION 2.13. *Let $\mathcal{A} \in \mathbb{R}^d$. Assume that $g:\mathbb{R}^d \to \mathbb{R}$ admits a decomposition $g = g_1 + g_2$ with $g_1$ a continuous function such that $\forall M > 0$,* $\mathbb{E}(\sup_{|\theta|\leq M}|g_1(G+\theta)|)<+\infty$ *and* $g_2\in\mathcal{V}_{\mathcal{A}}$. *Then, for any sequence* $(\vartheta_n)_n$ *of real-valued random variables converging a.s. to some deterministic limit* $\vartheta_\star\in\mathbb{R}$, $M_n(\mathcal{A}\vartheta_n,g)$ *converges a.s. to* $\mathbb{E}(g(G))$.

PROOF. By Proposition 2.6, it is enough to deal with the situation where $g=g_\uparrow+g_\downarrow$, with $g_\uparrow$ (resp., $g_\downarrow$) being an $\mathcal{A}$-nondecreasing (resp., $\mathcal{A}$-nonincreasing) function satisfying (2.1). One has $g=g_\uparrow\mathbf{1}_{\{g_\uparrow\geq 0\}}+g_\uparrow\mathbf{1}_{\{g_\uparrow<0\}}+g_\downarrow\mathbf{1}_{\{g_\downarrow\geq 0\}}+g_\downarrow\mathbf{1}_{\{g_\downarrow<0\}}$, where the functions $g_\uparrow\mathbf{1}_{\{g_\uparrow\geq 0\}}$ and $-g_\downarrow\mathbf{1}_{\{g_\downarrow<0\}}$ (resp., $g_\downarrow\mathbf{1}_{\{g_\downarrow\geq 0\}}$ and $-g_\uparrow\mathbf{1}_{\{g_\uparrow<0\}}$) are nonnegative, $\mathcal{A}$-nondecreasing (resp., $\mathcal{A}$-nonincreasing) and satisfy (2.1). As a consequence, it is enough to deal with the case where $g$ is nonnegative, $\mathcal{A}$-monotonic and satisfies (2.1). Choosing $\vartheta'\geq\vartheta$ when $g$ is $\mathcal{A}$-nondecreasing and $\vartheta\geq\vartheta'$ when $g$ is $\mathcal{A}$-nonincreasing, one has, for all $x\in\mathbb{R}^d$,

$$\begin{aligned}
&g(x+\mathcal{A}\vartheta')e^{-\mathcal{A}\vartheta'\cdot x-|\mathcal{A}\vartheta'|^2/2}-g(x+\mathcal{A}\vartheta)e^{-\mathcal{A}\vartheta\cdot x-|\mathcal{A}\vartheta|^2/2}\\
&\quad\geq (g(x+\mathcal{A}\vartheta')(e^{-\mathcal{A}\vartheta'\cdot x-|\mathcal{A}\vartheta'|^2/2}-e^{-\mathcal{A}\vartheta\cdot x-|\mathcal{A}\vartheta|^2/2}))\\
&\qquad\vee (g(x+\mathcal{A}\vartheta)(e^{-\mathcal{A}\vartheta'\cdot x-|\mathcal{A}\vartheta'|^2/2}-e^{-\mathcal{A}\vartheta\cdot x-|\mathcal{A}\vartheta|^2/2})).
\end{aligned}\tag{2.12}$$

From now on, we suppose that $g$ is nonnegative, $\mathcal{A}$-nondecreasing and satisfies (2.1): a symmetric argument applies to the nonincreasing case. Let $\varepsilon>0$, $\eta\in(0,1]$. For $m\in\mathbb{N}^*$,

$$\begin{aligned}
&\mathbb{P}(\exists n\geq m,|\mathbb{E}(g(G))-M_n(\mathcal{A}\vartheta_n,g)|\geq\varepsilon)\\
&\quad\leq\mathbb{P}\biggl(\exists n\geq m,|\mathbb{E}(g(G))-M_n(\mathcal{A}\vartheta_\star,g)|\geq\frac{\varepsilon}{2}\biggr)\\
&\qquad+\mathbb{P}(\exists n\geq m,|\vartheta_n-\vartheta_\star|>\eta)\\
&\qquad+\mathbb{P}\biggl(\forall n\geq m,|\vartheta_n-\vartheta_\star|\leq\eta\text{ and}\\
&\qquad\qquad\exists n\geq m,|M_n(\mathcal{A}\vartheta_\star,g)-M_n(\mathcal{A}\vartheta_n,g)|\geq\frac{\varepsilon}{2}\biggr).
\end{aligned}$$

By the strong law of large numbers and the a.s. convergence of $\vartheta_n$ to $\vartheta_\star$, the first two terms on the right-hand side both converge to 0 as $m\to+\infty$. Let us check that the third one also converges to 0. Let $M=|\vartheta_\star|+1$, $K>0$. For $-M\leq\vartheta\leq\vartheta'\leq M$, one has, using (2.12) for the first inequality, then (2.1) and (2.8), that

$$\begin{aligned}
&M_n(\mathcal{A}\vartheta',g)-M_n(\mathcal{A}\vartheta,g)\\
&\quad\geq-\frac{1}{n}\sum_{i=1}^{n}(|g(G_i+\mathcal{A}\vartheta')|\wedge|g(G_i+\mathcal{A}\vartheta)|)\\
&\qquad\qquad\times|e^{-\mathcal{A}\vartheta'\cdot G_i-|\mathcal{A}\vartheta'|^2/2}-e^{-\mathcal{A}\vartheta\cdot G_i-|\mathcal{A}\vartheta|^2/2}|\mathbf{1}_{\{|G_i|\leq K\}}
\end{aligned}$$

$$- \frac{1}{n} \sum_{i=1}^{n} (|g(G_i + \mathcal{A}\vartheta')| e^{-\mathcal{A}\vartheta' \cdot G_i - |\mathcal{A}\vartheta'|^2/2}$$

$$+ |g(G_i + \mathcal{A}\vartheta)| e^{-\mathcal{A}\vartheta \cdot G_i - |\mathcal{A}\vartheta|^2/2}) \mathbf{1}_{\{|G_i| > K\}}$$

$$\geq -\gamma_K(\vartheta' - \vartheta) - \frac{C}{n} \sum_{i=1}^{n} e^{c|G_i|^\beta + M|\mathcal{A}||G_i|} \mathbf{1}_{\{|G_i| > K\}},$$

where $\gamma_K = \lambda|\mathcal{A}| e^{c(K^\beta + (M|\mathcal{A}|)^\beta)} (M|\mathcal{A}| + K) e^{MK|\mathcal{A}|}$ and $C = 2\lambda e^{c(M|\mathcal{A}|)^\beta}$. When $|\vartheta_n - \vartheta_\star| \leq \eta$, choosing $\vartheta = \vartheta_n$ and $\vartheta' = \vartheta_\star + \eta$ then $\vartheta = \vartheta_\star - \eta$ and $\vartheta' = \vartheta_n$, one deduces that $M_n(\mathcal{A}\vartheta_\star, g) - M_n(\mathcal{A}\vartheta_n, g)$ is bounded from below and above, respectively, by

$$M_n(\mathcal{A}\vartheta_\star, g) - M_n(\mathcal{A}(\vartheta_\star + \eta), g) - \gamma_K(\vartheta_\star + \eta - \vartheta_n)$$
$$- \frac{C}{n} \sum_{i=1}^{n} e^{c|G_i|^\beta + M|\mathcal{A}||G_i|} \mathbf{1}_{\{|G_i| > K\}}$$

and

$$M_n(\mathcal{A}\vartheta_\star, g) - M_n(\mathcal{A}(\vartheta_\star - \eta), g) + \gamma_K(\vartheta_n + \eta - \vartheta_\star)$$
$$+ \frac{C}{n} \sum_{i=1}^{n} e^{c|G_i|^\beta + M|\mathcal{A}||G_i|} \mathbf{1}_{\{|G_i| > K\}}.$$

Choosing $K$ such that $\mathbb{E}(e^{c|G_i|^\beta + M|\mathcal{A}||G_i|} \mathbf{1}_{\{|G_i| > K\}}) \leq \frac{\varepsilon}{8C}$ and then $\eta$ such that $2\gamma_K \eta \leq \frac{\varepsilon}{8}$, we deduce that

$$\mathbb{P}\left(\forall n \geq m, |\vartheta_n - \vartheta_\star| \leq \eta \text{ and } \exists n \geq m, |M_n(\mathcal{A}\vartheta_\star, h) - M_n(\mathcal{A}\vartheta_n, h)| \geq \frac{\varepsilon}{2}\right)$$
$$\leq \mathbb{P}\left(\exists n \geq m, \frac{1}{n} \sum_{i=1}^{n} e^{c|G_i|^\beta + M|\mathcal{A}||G_i|} \mathbf{1}_{\{|G_i| > K\}} \geq \frac{\varepsilon}{4C}\right)$$
$$+ \mathbb{P}\left(\exists n \geq m, M_n(\mathcal{A}\vartheta_\star, g) - M_n(\mathcal{A}(\vartheta_\star + \eta), g) \leq -\frac{\varepsilon}{8}\right)$$
$$+ \mathbb{P}\left(\exists n \geq m, M_n(\mathcal{A}\vartheta_\star, g) - M_n(\mathcal{A}(\vartheta_\star - \eta), g) \geq \frac{\varepsilon}{8}\right).$$

By the strong law of large numbers $M_n(\mathcal{A}\vartheta_\star, g) - M_n(\mathcal{A}(\vartheta_\star + \eta), g)$ and $M_n(\vartheta_\star, g) - M_n(\mathcal{A}(\vartheta_\star - \eta), g)$ both converge a.s. to 0 and $\frac{1}{n} \sum_{i=1}^{n} e^{c|G_i|^\beta + M|\mathcal{A}||G_i|} \times \mathbf{1}_{\{|G_i| > K\}}$ to some limit not greater than $\frac{\varepsilon}{8C}$. One concludes that each term on the right-hand side converges to 0 as $m \to \infty$. $\square$

PROPOSITION 2.14. *Assume that* $g : \mathbb{R}^d \to \mathbb{R}$ *is such that* $\mathbb{E}(g^2(G + \theta_\star^{f,A}) \times e^{-2\theta_\star^{f,A} \cdot G}) < +\infty$ *and admits a decomposition* $g = g_1 + g_2 + \mathbf{1}_{\{d'=1\}} g_3$, *with* $g_1$ *of class* $C^1$ *satisfying* (2.2), $g_2 \in \mathcal{H}_\alpha$ *for* $\alpha \in (\frac{\sqrt{d'^2+8d'}-d'}{4}, 1]$ *and* $g_3 \in \mathcal{V}_A$. *Then, under* (0.1) *and* (1.2),

$$\sqrt{n}(M_n(\theta_n^{f,A}, g) - \mathbb{E}(g(G))) \xrightarrow{\mathcal{L}} \mathcal{N}_1(0, \mathrm{Var}(g(G + \theta_\star^{f,A})e^{-\theta_\star^{f,A}\cdot G - |\theta_\star^{f,A}|^2/2})).$$

As in Proposition 2.7, this statement is proved by combining the usual central limit theorem governing the convergence of $\sqrt{n}(M_n(\theta_\star^{f,A}, g) - \mathbb{E}(g(G)))$, Lemma 2.8, Proposition 2.9, the decomposition of functions in $\mathcal{V}_A$ given at the beginning of the proof of Proposition 2.13 and the next result.

PROPOSITION 2.15. *Let* $\mathcal{A} \in \mathbb{R}^d$ *and* $g : \mathbb{R}^d \to \mathbb{R}$ *be an* $\mathcal{A}$*-monotonic function with constant sign satisfying* (2.1),

$$\forall \delta > 1/4, \forall \vartheta_0 \in \mathbb{R}, \qquad \sup_{\vartheta \in [\vartheta_0 \pm 1/n^\delta]} \sqrt{n}|M_n(\mathcal{A}\vartheta, g) - M_n(\mathcal{A}\vartheta_0, g)| \xrightarrow{\mathrm{Pr}} 0.$$

REMARK 2.16. Assume that $d' = 1$. Let $g \in \mathcal{V}_A$, and $(\nu_n)_n$ be a deterministic integer-valued sequence such that

$$\exists \lambda > 0, \exists \gamma > \tfrac{1}{2}, \forall n \in \mathbb{N}^*, \qquad \nu_n \geq \lambda n^\gamma.$$

Combining Propositions 1.2 and 2.15, one obtains that under (0.1) and (1.2), $\sqrt{n}(M_n(\theta_{\nu_n}^{f,A}, g) - \mathbb{E}(g(G))) \xrightarrow{\mathcal{L}} \mathcal{N}_1(0, \mathrm{Var}(g(G + \theta_\star^{f,A})e^{-\theta_\star^{f,A}\cdot G - |\theta_\star^{f,A}|^2/2}))$.

PROOF OF PROPOSITION 2.15. Up to a multiplication by $-1$, we may assume that $g$ is nonnegative. Moreover, we only deal with the case where $g$ is $\mathcal{A}$-nondecreasing, the nonincreasing case being obtained by a symmetric argument. By (2.12), for $\vartheta' < \vartheta''$ and $\vartheta \in [\vartheta', \vartheta'']$,

$$\begin{aligned} M_n(\mathcal{A}\vartheta', g) &- \frac{1}{n}\sum_{i=1}^n |g(G_i + \mathcal{A}\vartheta)||e^{-\mathcal{A}\vartheta\cdot G_i - |\mathcal{A}\vartheta|^2/2} - e^{-\mathcal{A}\vartheta'\cdot G_i - |\mathcal{A}\vartheta'|^2/2}| \\ &\leq M_n(\mathcal{A}\vartheta, g) \\ &\leq M_n(\mathcal{A}\vartheta'', h) \\ &\quad + \frac{1}{n}\sum_{i=1}^n |g(G_i + \mathcal{A}\vartheta)||e^{-\mathcal{A}\vartheta\cdot G_i - |\mathcal{A}\vartheta|^2/2} - e^{-\mathcal{A}\vartheta''\cdot G_i - |\mathcal{A}\vartheta''|^2/2}|. \end{aligned}$$

With (2.1) and (2.8), one deduces that if $-M \leq \vartheta' \leq \vartheta'' \leq M$, then

$$\begin{aligned} \sup_{\vartheta \in [\vartheta', \vartheta'']} &|M_n(\mathcal{A}\vartheta, g) - M_n(\mathcal{A}\vartheta_0, g)| \\ &\leq \max(|M_n(\mathcal{A}\vartheta', g) - M_n(\mathcal{A}\vartheta_0, g)|, |M_n(\mathcal{A}\vartheta'', g) - M_n(\mathcal{A}\vartheta_0, g)|) \\ &\quad + \frac{C(\vartheta'' - \vartheta')}{n}\sum_{i=1}^n e^{c|G_i|^\beta + M|\mathcal{A}||G_i|}(M|\mathcal{A}| + |G_i|). \end{aligned} \tag{2.13}$$

Let $\nu = \frac{\beta\vee 1}{2}$ and $M = |\vartheta_0| + 1$. When $\max_{1\le i\le n}|G_i| \le \sqrt{2d\log n}$, the second term on the right-hand side is smaller than $\gamma e^{\gamma(\log n)^\nu}(\vartheta'' - \vartheta')$, where the constant $\gamma$ does not depend on $n$. Let $\varepsilon > 0$. We set $K = \lceil 2\gamma n^{1/2-\delta} e^{\gamma(\log n)^\nu}/\varepsilon\rceil$ and $\vartheta_k = \vartheta_0 + k\varepsilon/2\gamma e^{\gamma(\log n)^\nu}$ for $k \in \{-K,\dots,K\}$. Applying (2.13) with $\vartheta' = \vartheta_k$ and $\vartheta'' = \vartheta_{k+1}$, one obtains that when $\max_{1\le i\le n}|G_i| \le \sqrt{2d\log n}$,

$$\begin{aligned}\sup_{\vartheta\in[\vartheta_k,\vartheta_{k+1}]} &|M_n(\mathcal{A}\vartheta, g) - M_n(\mathcal{A}\vartheta_0, g)|\\ &\le \frac{\varepsilon}{2\sqrt{n}} + \max(|M_n(\mathcal{A}\vartheta_k, g) - M_n(\mathcal{A}\vartheta_0, g)|,\\ &\qquad\qquad |M_n(\mathcal{A}\vartheta_{k+1}, g) - M_n(\mathcal{A}\vartheta_0, g)|).\end{aligned}$$

Therefore,

$$\begin{aligned}\mathbb{P}\Big(&\sup_{\vartheta\in[\vartheta_0\pm 1/n^\delta]} |M_n(\mathcal{A}\vartheta, g) - M_n(\mathcal{A}\vartheta_0, g)| \ge \frac{\varepsilon}{\sqrt{n}}\Big)\\ &\le \mathbb{P}\Big(\max_{1\le i\le n}|G_i| > \sqrt{2d\log n}\Big)\\ &\quad + \mathbb{P}\Big(\max_{1\le i\le n}|G_i| \le \sqrt{2d\log n},\\ &\qquad\quad \max_{|k|\le K}|M_n(\mathcal{A}\vartheta_k, g) - M_n(\mathcal{A}\vartheta_0, g)| \ge \frac{\varepsilon}{2\sqrt{n}}\Big).\end{aligned}$$

By (2.9), the first term on the right-hand side tends to 0 as $n\to\infty$. Reasoning as we did at the end of the proof of Proposition 2.9, with the next lemma replacing Lemma 2.12, we conclude that the second term also tends to 0. $\square$

Lemma 2.17. *When $\mathcal{A}\in\mathbb{R}^d$ and $g:\mathbb{R}^d\to\mathbb{R}$ is a $\mathcal{A}$-monotonic function with constant sign satisfying (2.1), we have*

$$\begin{aligned}&\forall M > 0, \exists C > 0, \forall\vartheta,\vartheta' \in [-M, M], \forall n\in\mathbb{N}^*\\ &\mathbb{E}((M_n(\mathcal{A}\vartheta, g) - M_n(\mathcal{A}\vartheta', g))^2) \le \frac{C|\vartheta - \vartheta'|}{n}.\end{aligned}$$

Proof. Choosing $\vartheta' \ge \vartheta$ if $g$ is nonnegative and $\mathcal{A}$-nondecreasing, or nonpositive and $\mathcal{A}$-nonincreasing, and $\vartheta \ge \vartheta'$ otherwise, one has

$$\begin{aligned}\mathbb{E}((&g(G+\mathcal{A}\vartheta)e^{-\mathcal{A}\vartheta\cdot G - |\mathcal{A}\vartheta|^2/2} - g(G+\mathcal{A}\vartheta')e^{-\mathcal{A}\vartheta'\cdot G - |\mathcal{A}\vartheta'|^2/2})^2)\\ &= \mathbb{E}(g^2(G)e^{-\mathcal{A}\vartheta\cdot G + |\mathcal{A}\vartheta|^2/2}) + \mathbb{E}(g^2(G)e^{-\mathcal{A}\vartheta'\cdot G + |\mathcal{A}\vartheta'|^2/2})\\ &\quad - 2\mathbb{E}(g(G)g(G+\mathcal{A}(\vartheta'-\vartheta))e^{-\mathcal{A}\vartheta'\cdot G + \mathcal{A}\vartheta\cdot\mathcal{A}\vartheta' - |\mathcal{A}\vartheta'|^2/2})\end{aligned}$$

$$\le \mathbb{E}(g^2(G)(e^{-\mathcal{A}\vartheta\cdot G+|\mathcal{A}\vartheta|^2/2} + e^{-\mathcal{A}\vartheta'\cdot G+|\mathcal{A}\vartheta'|^2/2} - 2e^{-\mathcal{A}\vartheta'\cdot G+\mathcal{A}\vartheta\cdot\mathcal{A}\vartheta'-|\mathcal{A}\vartheta'|^2/2})).$$

The conclusion is then a consequence of the following inequality: for $\theta, \theta' \in \mathbb{R}^d$ with $|\theta| \vee |\theta'| \le |\mathcal{A}|M$,

$$\begin{aligned}
&\mathbb{E}(g^2(G)(e^{-\theta\cdot G+|\theta|^2/2} + e^{-\theta'\cdot G+|\theta'|^2/2} - 2e^{-\theta'\cdot G+\theta\cdot\theta'-|\theta'|^2/2})) \\
&\quad\le C\mathbb{E}(e^{|G|^2/4}(|e^{-\theta\cdot G+|\theta|^2/2} - e^{-\theta'\cdot G+|\theta'|^2/2}| \\
&\qquad\qquad + 2e^{-\theta'\cdot G-|\theta'|^2/2}|e^{|\theta'|^2} - e^{\theta\cdot\theta'}|)) \\
&\quad\le C\biggl(|\theta-\theta'|\int_0^1 e^{3|\vartheta(t)|^2/2}\int_{\mathbb{R}^d}|\vartheta(t)-x|e^{-|x+2\vartheta(t)|^2/4}\,dx\,dt \\
&\qquad\qquad + 2e^{|\theta'|^2/2}|e^{|\theta'|^2} - e^{\theta\cdot\theta'}|\int_{\mathbb{R}^d} e^{-|x+2\theta'|^2/4}\,dx\biggr) \\
&\quad\le C|\theta-\theta'|.
\end{aligned}$$

□

**3. Practical implementation and applications.** Option pricing in local or stochastic volatility models eventually boils down to the computation of an expectation $\mathbb{E}(f(G))$, where $G$ is a $d$-dimensional standard normal random vector. In a financial context, there is no restriction in assuming that the payoff function $f$ satisfies both (0.1) and (0.2). In most cases, this expectation will be computed using Monte Carlo simulations because closed formulas are barely available. The question of reducing the variance arises quite naturally in this context. Relying on equation (0.3), we have chosen the importance sampling point of view to tackle the delicate problem of variance reduction. Practitioners' desires with variance reduction is to have an automatic toolbox at hand, which is precisely what we are devising here. As explained in the Introduction, we advise to compute the minimizer $\vartheta_n^{f,A}$ of $v_n^{f,A}$ and then to use this value in a Monte Carlo procedure, as described in Algorithm 1. Note that the same samples are used to compute $\vartheta_n^{f,A}$ and the Monte Carlo estimator $M_n(A\vartheta_n^{f,A}, f)$. Even though the terms involved in $M_n^f(A\vartheta_n^{f,A}, f)$ are not independent, according to Corollary 2.4, it is as easy to construct confidence intervals, as for a crude Monte Carlo computation.

REMARK 3.1. In the name ("Reduced Robust Importance Sampling") of Algorithm 1, the term "Reduced" emphasizes that the optimal importance sampling parameter is searched for in a subspace of the set of all parameters. When the matrix $A = I_d$, the algorithm is simply denoted RIS because there is no longer any dimension reduction.

In this section, we first explain how $\vartheta_n^{f,A}$ can be computed using Newton's optimization procedure. We then illustrate the efficiency of this robust

**Algorithm 1** Reduced Robust Importance Sampling (RRIS)

1. Generate $G_1, \dots, G_n$, $n$ i.i.d. samples following the law of $G$.
2. Compute the minimizer $\vartheta_n^{f,A}$ of

$$v_n^{f,A}(\vartheta) = \frac{1}{n}\sum_{i=1}^{n} f^2(G_i) e^{-A\vartheta\cdot G_i + |A\vartheta|^2/2}.$$

3. Compute the expectation $\mathbb{E}(f(G))$ by Monte Carlo

$$M_n(A\vartheta_n^{f,A}, f) = \frac{1}{n}\sum_{i=1}^{n} f(G_i + A\vartheta_n^{f,A}) e^{-A\vartheta_n^{f,A}\cdot G_i - |A\vartheta_n^{f,A}|^2/2}.$$

variance reduction technique, both in the multidimensional Black–Scholes framework and in more general local volatility frameworks.

3.1. *Solving the minimization problem.* We already know, from Proposition 1.1, that the function $v_n^{f,A}$ is strongly convex and infinitely continuously differentiable. Hence, we can approximate $\vartheta_n^{f,A}$ using Newton's algorithm, for instance. The Hessian matrix $\nabla^2_\vartheta v_n^{f,A}(\vartheta)$ can be written as the sum of a scalar matrix and a positive semidefinite matrix. Hence, it is quite obvious that the smallest eigenvalue of $\nabla^2_\vartheta v_n^{f,A}(\vartheta)$ is larger than the smallest eigenvalue of $A^*A$ times $\frac{1}{n}\sum_{i=1}^n f^2(G_i)e^{-A\vartheta\cdot G_i+|A\vartheta|^2/2}$. This last term can be arbitrarily small, depending on the function $f$. Therefore a straightforward application of Newton's algorithm can be particularly inefficient in some cases. It would be much better to have an alternative representation of $\vartheta_n^{f,A}$ as the minimizer of a function, the smallest eigenvalue of whose Hessian matrix does not depend on $f$. We advise to rewrite $\nabla_\vartheta v_n^{f,A}(\vartheta)$ as

$$\nabla_\vartheta v_n^{f,A}(\vartheta) = A^*A\vartheta \frac{1}{n}\sum_{i=1}^{n} f^2(G_i) e^{-A\vartheta\cdot G_i+|A\vartheta|^2/2}$$
$$- \frac{1}{n}\sum_{i=1}^{n} A^* G f^2(G_i) e^{-A\vartheta\cdot G_i+|A\vartheta|^2/2}.$$

Hence, $\vartheta_n^{f,A}$ can be seen as the root of

$$\nabla_\vartheta u_n^{f,A}(\vartheta) = A^*A\vartheta - \frac{\sum_{i=1}^n A^* G_i f^2(G_i) e^{-A\vartheta\cdot G_i}}{\sum_{i=1}^n f^2(G_i) e^{-A\vartheta\cdot G_i}},$$

with $u_n^{f,A}(\vartheta) = \frac{|A\vartheta|^2}{2} + \log(\sum_{i=1}^n f^2(G_i)e^{-A\vartheta\cdot G_i})$. The Hessian matrix of $u_n^{f,A}$ is given by

$$\nabla^2_\vartheta u_n^{f,A}(\vartheta) = A^*A + \frac{\sum_{i=1}^n A^* G_i G_i^* A f^2(G_i) e^{-A\vartheta\cdot G_i}}{\sum_{i=1}^n f^2(G_i) e^{-A\vartheta\cdot G_i}}$$

$$- \frac{(\sum_{i=1}^n A^* G_i f^2(G_i) e^{-A\vartheta\cdot G_i})(\sum_{i=1}^n A^* G_i f^2(G_i) e^{-A\vartheta\cdot G_i})^*}{(\sum_{i=1}^n f^2(G_i) e^{-A\vartheta\cdot G_i})^2}.$$

Using the Cauchy–Schwarz inequality, it is clear that $\nabla_\vartheta^2 u_n^{f,A}(\vartheta) - A^*A$ is a positive semidefinite matrix. Hence, the smallest eigenvalue of $\nabla_\vartheta^2 u_n^{f,A}(\vartheta)$ is always larger than the smallest one of $A^*A$, whatever the values taken by $f$ are. This advocates the use of $u_n^{f,A}$, rather than $v_n^{f,A}$, to compute $\vartheta_n^{f,A}$.

Using this new expression, we implement Algorithm 2 to construct an approximation $x_n^k$ of $\vartheta_n^{f,A}$. Since $u_n^{f,A}$ is strongly convex, for any fixed $n$, $x_n^k$ converges to $\vartheta_n^{f,A}$ when $k$ goes to infinity. The direction of descent $d_n^k$ at step $k$ should be computed as the solution of a linear system. There is no point in computing the inverse of $\nabla_\vartheta^2 u_n^{f,A}(x_n^k)$, which would be computationally much more expensive.

Remarks on the implementation. From a practical point of view, $\varepsilon$ should be chosen reasonably small, $\varepsilon \approx 10^{-6}$. This algorithm converges very quickly and, in most cases, less than five iterations are enough to get a very accurate estimate of $\vartheta_n^{f,A}$, actually within the $\varepsilon$-error. Since the points at which the payoff function $f$ is evaluated remain constant through the iterations of Newton's algorithm, the values $f^2(G_i)$ for $i = 1, \ldots, n$ should be computed before starting the optimization algorithm, something which considerably speeds up the whole process. The Hessian matrix of our problem is easily tractable, so there is no point in using quasi-Newton's methods.

3.2. *Numerical examples.* In this subsection, we present numerical results obtained by combining Algorithms 1 and 2 for different pricing problems. For each example, we have computed the reference price using a crude Monte Carlo estimator with a huge number of samples, such that the width of the 95% confidence interval is $10^{-3}$. In the columns "price MC," "price RIS" and "price RRIS" of the tables, we give, respectively, the crude Monte Carlo estimator, the RIS estimator and the RRIS estimator of the price computed with the same, smaller, number $n$ of samples. The variances for the RRIS algorithm (resp., the RIS algorithm) given in the

**Algorithm 2** Newton's algorithm

Choose an initial value $x_n^0 \in \mathbb{R}^d$.
$k = 1$
**while** $|\nabla_\vartheta u_n^{f,A}(x_n^k)| > \varepsilon$ **do**
 1. Compute $d_n^k$ such that $(\nabla_\vartheta^2 u_n^{f,A}(x_n^k)) d_n^k = -\nabla_\vartheta u_n^{f,A}(x_n^k)$.
 2. $x_n^{k+1} = x_n^k + d_n^k$, $k = k + 1$.
**end while**

tables below are computed along a single run of the algorithm using the estimator $v_n^{f,A}(\vartheta_n^{f,A}) - M_n^2(\theta_n^{f,A}, f)$ [resp., $v_n^f(\theta_n^f) - M_n^2(\theta_n^f, f)$] which converges almost surely to $v^f(\theta_\star^f) - \mathbb{E}^2(f(G))$ under the assumptions of Theorem 2.2. The variances of the crude Monte Carlo methods (denoted 'Var MC' in the tables) are estimated by $\frac{1}{n}\sum_{i=1}^n f^2(G_i) - (\frac{1}{n}\sum_{i=1}^n f(G_i))^2$. All of the histograms presented hereafter are centered around their empirical means and renormalized by the empirical variances. When no further indications are given, the matrix $A$ is chosen as the identity, which implies that $d = d'$ and $\theta_n^{f,A} = \vartheta_n^{f,A}$.

3.2.1. *Black–Scholes framework.* First, we consider an $I$-dimensional Black–Scholes model in which the dynamics under the risk-neutral measure of each asset $S^i$ is supposed to be given by

$$dS_t^i = S_t^i(r\,dt + \sigma^i\,dW_t^i), \qquad S_0 = (S_0^1, \dots, S_0^I),$$

where $W = (W^1, \dots, W^I)$. Each component $W^i$ is a standard Brownian motion. For the numerical experiments, the covariance structure of $W$ will be assumed to be given by $\langle W^i, W^j\rangle_t = \rho t \mathbf{1}_{\{i\neq j\}} + t\mathbf{1}_{\{i=j\}}$. We suppose that $\rho \in (-\frac{1}{I-1}, 1)$, which ensures that the matrix $C = (\rho\mathbf{1}_{\{i\neq j\}} + \mathbf{1}_{\{i=j\}})_{1\le i,j\le I}$ is positive definite. Let $L$ denote the lower-triangular matrix involved in the Cholesky decomposition $C = LL^*$. To simulate $W$ on the time-grid $0 < t_1 < t_2 < \cdots < t_N$, we need $d = I \times N$ independent standard normal variables and set

$$\begin{pmatrix} W_{t_1} \\ W_{t_2} \\ \vdots \\ W_{t_{N-1}} \\ W_{t_N} \end{pmatrix} = \begin{pmatrix} \sqrt{t_1}L & 0 & 0 & \dots & 0 \\ \sqrt{t_1}L & \sqrt{t_2 - t_1}L & 0 & \dots & 0 \\ \vdots & \ddots & \ddots & \ddots & \vdots \\ \vdots & \ddots & \ddots & \sqrt{t_{N-1} - t_{N-2}}L & 0 \\ \sqrt{t_1}L & \sqrt{t_2 - t_1}L & \dots & \sqrt{t_{N-1} - t_{N-2}}L & \sqrt{t_N - t_{N-1}}L \end{pmatrix} G,$$

where $G$ is a normal random vector in $\mathbb{R}^{I\times N}$. The vector $(\sigma^1, \dots, \sigma^d)$ is the vector of volatilities and $r > 0$ is the instantaneous interest rate. We will denote the maturity time by $T$.

*Basket option.* We consider options with payoffs of the form $(\sum_{i=1}^d \omega^i S_T^i - K)_+$, where $(\omega^1, \dots, \omega^d)$ is a vector of algebraic weights. The strike value $K$ can be taken to be negative, to deal with put-like options. All of these payoffs belong to $\mathcal{H}_1$, so Theorem 2.3 applies, as Figures 1 and 2 illustrate. These histograms have been obtained with 5000 independent runs of the RIS algorithm. The case of such basket options is definitely a crucial issue because there is no closed formula as soon as $d > 2$, and the variance of a crude Monte Carlo approach can be dramatically large. We can see in the examples of the basket options treated in Table 1 that the Robust Importance Sampling method does reduce the variance by at least 10. The results are obtained within 4.5 CPU seconds, compared to the 1.5 CPU seconds needed for the crude Monte Carlo computation. The same number of samples are used in both methods, which brings an overall gain of 3.3 in favor of the RIS algorithm. In the case $\rho = 0.2$ and $K = 50$, which is the option used for the histograms, the empirical variance is 1.76, whereas the on-line estimated variance is 1.74. This illustrates the conclusion of Corollary 2.4. The improvement brought by the RIS algorithm is very encouraging, not only because it definitely reduces the variance, but, above all, because it is fully automatic. Unlike most adaptive importance sampling strategies developed so far and, in particular, the ones based on stochastic approximations, the one we propose here does not require any parameter tuning.

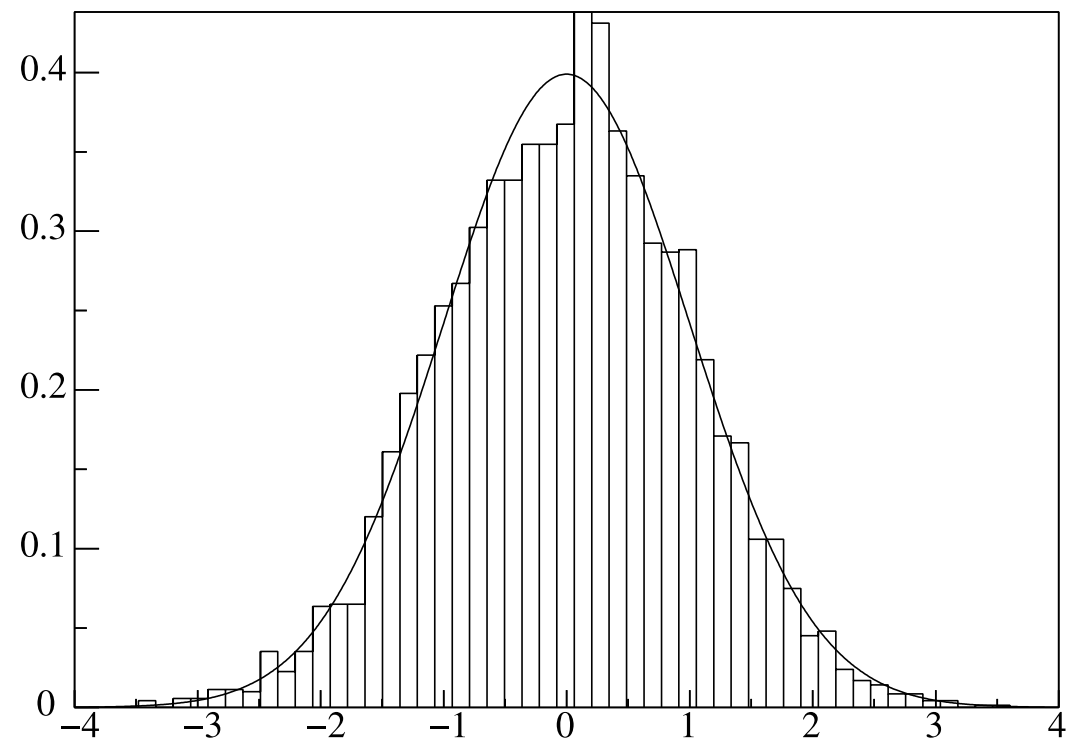

FIG. 1. *Limiting distribution of* $\theta_n^f$ *for the option of Table 1 with* $\rho = 0.2$ *and* $K = 50$.

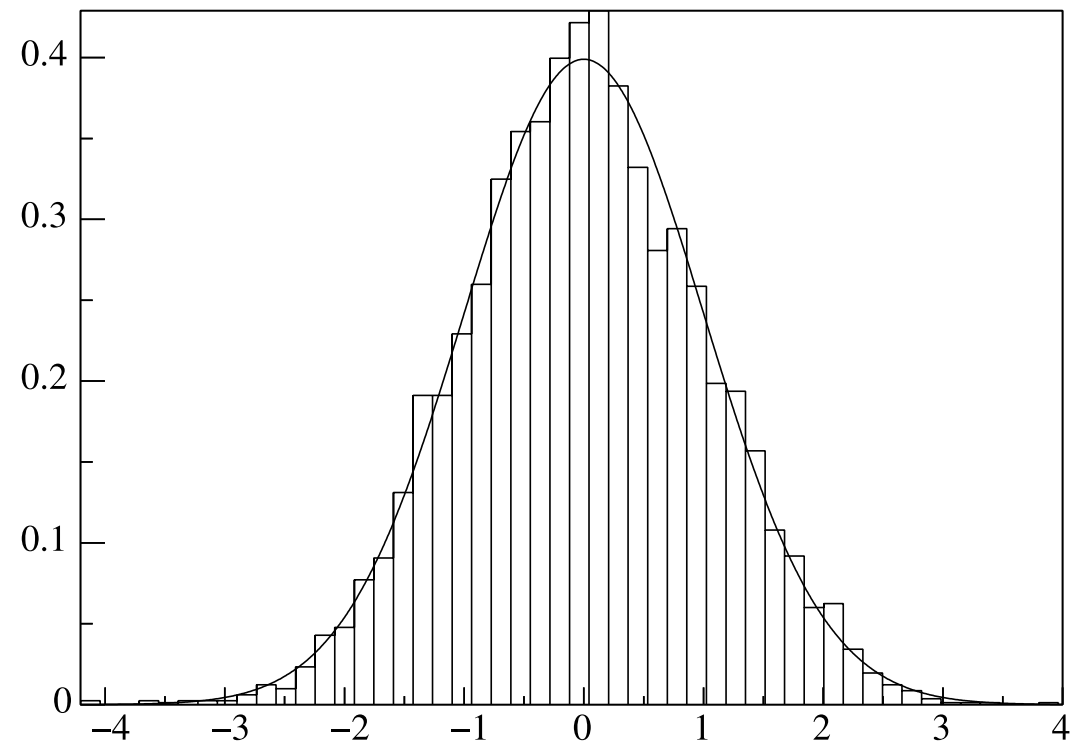

FIG. 2. *Limiting distribution of* $M_n(\theta_n^f, f)$ *(RIS) for the option of Table 1 with* $\rho = 0.2$ *and* $K = 50$.

We have also tested our algorithm on a 10-dimensional exchange option with randomly chosen spots and volatilities. The numerical results of Table 2 show that the RIS algorithm performs well for a wide variety of basket options. In any case, the variance is divided by at least 7, whereas it increases twice the CPU time. This leads to an overall gain of 3.5 in the worst case.

*One-dimensional digital option.* We consider an option with payoff $\mathbf{1}_{\{S_T>L\}}$, where $L>0$. We choose $T=1$, $S_0=100$, $\sigma=0.2$, $r=0.05$ and $L=140$. We fix the number of samples at 100,000. A crude Monte Carlo computation gives a price of 0.05952 with a variance of 0.053, whereas the exact price is 0.05968. On each run of the algorithm, we can compute the on-line estimator of the variance and use it to construct a confidence interval. We have run the RIS algorithm 100,000 times independently and, on each run, we have constructed the confidence interval of level 95% using the on-line estimated variance $v_n^{f,A}(\vartheta_n^{f,A}) - M_n^2(\theta_n^{f,A}, f)$. The true price falls outside the confidence interval in 5104 cases out of 100,000, which gives a level of 94.9%. This little experiment illustrates how Corollary 2.4 can be used to construct confidence intervals.

Table 1

*Basket option in dimension $d=40$ with $r=0.05$, $T=1$, $S_0^i=50$, $\sigma^i=0.2$, $\omega^i=\frac{1}{d}$ for all $i=1,\dots,d$ and $n=10{,}000$*

| $\rho$ | $K$ | Price | Price MC | Variance MC | Price RIS | Variance RIS |
|---|---|---|---|---|---|---|
| 0.1 | 45 | 7.210 | 7.216 | 12.12 | 7.209 | 1.04 |
| | 55 | 0.561 | 0.567 | 1.90 | 0.559 | 0.14 |
| 0.2 | 50 | 3.298 | 3.304 | 13.56 | 3.296 | 1.74 |
| 0.5 | 45 | 7.662 | 7.678 | 42.2 | 7.650 | 5.06 |
| | 55 | 1.906 | 1.879 | 14.46 | 1.906 | 1.25 |
| 0.9 | 45 | 8.215 | 8.154 | 69.47 | 8.211 | 7.89 |
| | 55 | 2.823 | 2.823 | 30.08 | 2.819 | 2.58 |

Table 2

*Basket option in dimension $d=10$ with $r=0.05$, $T=1$, $K=0$, $\rho=0.2$. The spots are chosen uniformly in $[70,130]$ and the volatilities in $[0.1,0.3]$. $\omega^i=\frac{1}{d}$ for $i=1,\dots,d/2$ and $\omega^i=-\frac{1}{d}$ for $i=d/2+1,\dots,d$ and $n=100{,}000$*

| Price | Variance MC | Variance RIS |
|---|---|---|
| 3.58 | 21.66 | 2.97 |
| 0.129 | 0.511 | 0.016 |
| 7.4 | 34.04 | 5.02 |
| 1.08 | 5.24 | 0.52 |

*One dimensional barrier option.* This time, we only focus on one asset and we want to price a call option with a discrete barrier on this asset. A discrete barrier means that we only check if the asset has crossed the barrier at fixed dates $t_1, \dots, t_d = T$, usually one per month. We assume that the grid defined by $t_1, \dots, t_d$ is regular with step size $\delta t = T/d$. The payoff can be written as $(S_T - K)_+ \mathbf{1}_{\{\forall 1 \le i \le d, S_{t_i} \ge L\}}$ for a down-and-out call option with barrier $L$. The price of such an option can be written as $\mathbb{E}(f(G^1, \dots, G^d))$ with

$$f(x_1, \dots, x_d) = e^{-rT}(S_0 e^{(r-\sigma^2/2)T + \sigma\sqrt{\delta t}\sum_{j=1}^d x_j} - K)_+ \\ \times \mathbf{1}_{\{\forall 1 \le i \le d, S_0 e^{(r-\sigma^2/2)t_i + \sigma\sqrt{\delta t}\sum_{j=1}^i x_j} \ge L\}}.$$

In this particular case, if we consider the RIS algorithm developed before, the importance sampling parameter $\theta$ lies in $\mathbb{R}^d$. Hence, the optimization problem becomes harder to solve as the number of time steps increases.

One idea is to restrict the parameter $\theta$ to the subspace $\{A\vartheta : \vartheta \in \mathbb{R}\}$, where the vector $A$ is defined by $A = (\sqrt{t_1}, \dots, \sqrt{t_d - t_{d-1}})^*$. In this case, the optimal parameter is always real-valued $d' = 1$, whatever the number of time steps we consider. This alternative approach—referred to as RRIS (Reduced Robust Importance Sampling)—corresponds to adding a linear drift to the Brownian motion. These two approaches are compared in Table 3 for the case of a down-and-out call option and it turns out that the optimal variances obtained in both cases are very close to each other. When the underlying asset is of dimension one, the computation time gained by using the RRIS algorithm instead of the RIS one is not that important, but it will become a crucial issue for multidimensional barrier options. The efficiency of the two algorithms on the down-and-out call option is very impressive. As in the previous example, the variance is reduced by a factor between 8 and 11. The use of the RRIS algorithm compared to a crude Monte Carlo method doubles the computation time, which means that the gain is at least 4. Figures 3 and 4 illustrate the asymptotic behavior of the RIS algorithm. They have been obtained by running the RIS algorithm 5000 times independently. The histogram of Figure 3 represents the limiting distribution of the first component of $\theta_n^f$ computed with the RIS algorithm and fits the density of the standard normal distribution (plain line) rather well, which illustrates Proposition 1.2. Although the hypotheses of Theorems 2.2 and 2.3 are not satisfied for the payoff at hand in the RIS framework, Figure 4 shows that our estimator is still convergent and asymptotically normal. This numerical convergence is emphasized by the matching of the empirical variance of the histogram and the on-line variance computed on a single run of the RIS algorithm; for these two quantities, we find, respectively, 34.70 and 35.68. Since the payoff belongs to $\mathcal{V}_A$, the convergence and the asymptotic normality of the RRIS estimator are, in return, ensured by Theorems 2.2 and 2.3.

TABLE 3
*Down-and-out call option with* $\sigma = 0.2$, $r = 0.05$, $T = 2$, $S_0^1 = 100$, $K = 110$ *and* $n = 10{,}000$

| $L$ | Price | Price MC | Variance MC | Variance RIS | Price RRIS | Variance RRIS |
|---|---|---|---|---|---|---|
| 70 | 11.445 | 11.472 | 401.51 | 34.10 | 11.454 | 34.33 |
| 80 | 11.244 | 11.240 | 401.04 | 35.68 | 11.261 | 36.11 |
| 90 | 9.689 | 9.672 | 383.93 | 42.54 | 9.705 | 45.37 |
| 95 | 7.564 | 7.518 | 342.05 | 42.01 | 7.557 | 49.84 |

*Barrier basket option.* We consider basket options in dimension $I$ with a discrete barrier on each asset. For instance, if we consider a down-and-out call option, the payoff can be written as $(\sum_{i=1}^{I} \omega^i S_T^i - K)_+ \mathbf{1}_{\{\forall i \le I, \forall j \le N, S_{t_j}^i \ge L^i\}}$, where $\omega = (\omega^1, \dots, \omega^I)$ is a vector of positive weights, $L = (L^1, \dots, L^I)$ is the vector of barriers, $K > 0$ the strike value and $t_N = T$. Once again, we consider one time step per month, which means that for an option with maturity time $T = 2$ as in Table 4, the number of time steps is $N = 24$. From now on, we fix $I = 5$. Hence, in the RIS algorithm, the parameter $\theta$ is of dimension $d = 120$. Even though this is not that huge, it requires much more computational time, as the numerical experiments show. For the option of Table 4, a standard Monte Carlo computation takes 4.3 CPU seconds and the RRIS algorithm takes 8.7 CPU seconds, whereas the RIS algorithm needs 22.5 CPU seconds. The RIS algorithm is three times slower than the RRIS algorithm, in which the parameter $\theta$ lies in the subspace $\{A\vartheta : \vartheta \in \mathbb{R}^d\}$ of dimension $d' = I = 5$ with $A_{(j-1)I+i,i} = \sqrt{t_j - t_{j-1}}$ (convention $t_0 = 0$) for $j = 1, \dots, N$ and $i = 1, \dots, I$, all the other coefficients of $A$ being zero.

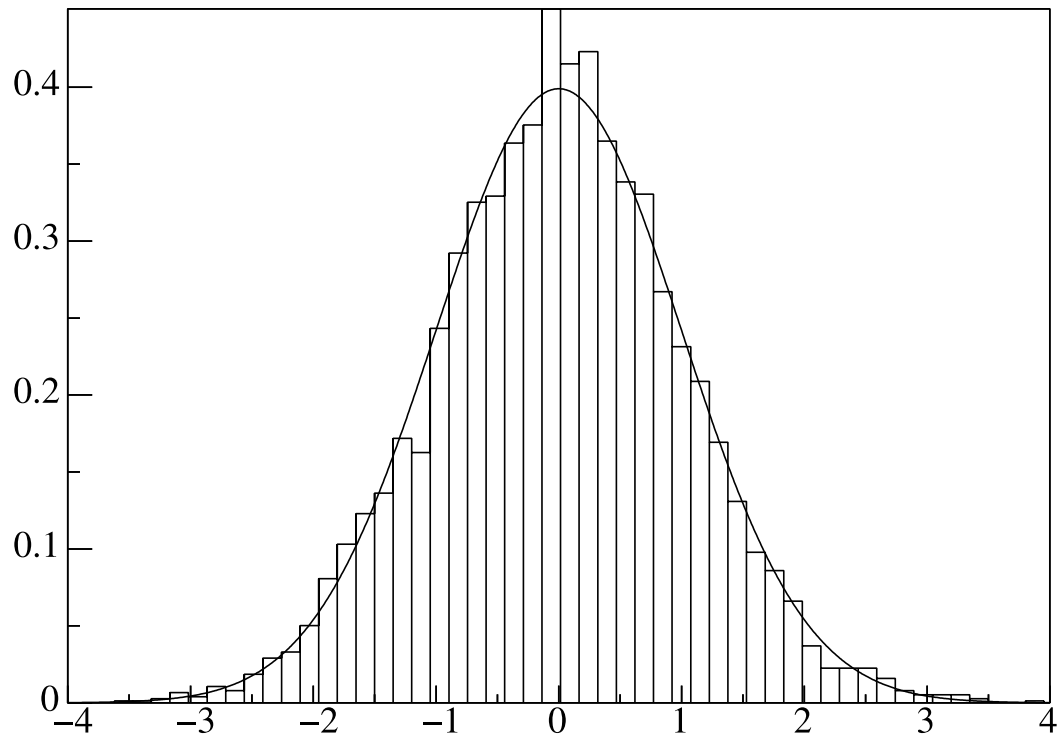

FIG. 3. *Limiting distribution of the first component of* $\theta_n^f$ *(RIS) for the option of Table 3 with* $L = 80$.

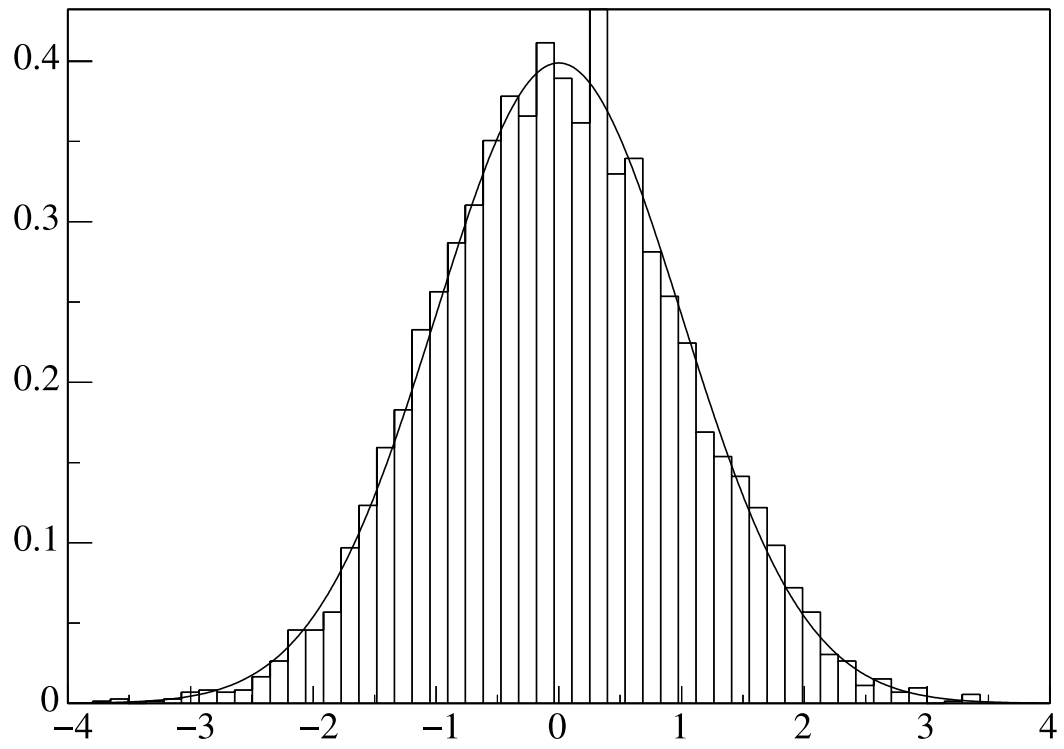

FIG. 4. *Limiting distribution of* $M_n(\theta_n^f, f)$ *(RIS) for the option of Table 3 with* $L = 80$.

Path-dependent basket options are a prime example of pricing problems in which the use of one importance sampling parameter per time step dramatically slows down the computation. Restricting the importance sampling parameter space to a subspace of dimension $d' = I = 5$, as in the RRIS algorithm, divides the computational time by 3, whereas the optimal variance of the RRIS algorithm is very close to that of the RIS algorithm. Hence, there is no point in using one importance sampling parameter per time step. The improvement factor in terms of variance provided by the RRIS algorithm varies between 10 and 20. Because the RRIS algorithm is twice as slow as a standard Monte Carlo computation, the overall gain factor varies between 5 and 10.

The payoff does not satisfy the assumptions of Theorems 2.2 and 2.3, neither in the RIS nor in the RRIS framework. Nevertheless, it seems rather clear from Figure 6 that the RRIS estimator is convergent and asymptotically normal. Besides, for $K = 50$, the variance computed on a single run of the RRIS algorithm perfectly matches the empirical variance of the histogram. Figure 5 illustrates the asymptotic normality of $\theta_n^{f,A}$ which is still ensures by Proposition 1.2 in this example. These histograms have been constructed from 100,000 independent runs of the RRIS algorithm.

TABLE 4

*Down-and-out call option in dimension* $I = 5$ *with* $\sigma = 0.2$, $S_0 = (50, 40, 60, 30, 20)$, $L = (40, 30, 45, 20, 10)$, $\rho = 0.3$, $r = 0.05$, $T = 2$, $\omega = (0.2, 0.2, 0.2, 0.2, 0.2)$ *and* $n = 100{,}000$

| $K$ | Price | Price MC | Var MC | Var RIS | Price RRIS | Var RRIS |
|---|---|---|---|---|---|---|
| 45 | 2.371 | 2.348 | 22.46 | 2.58 | 2.378 | 2.62 |
| 50 | 1.175 | 1.178 | 10.97 | 0.78 | 1.179 | 0.79 |
| 55 | 0.515 | 0.513 | 4.72 | 0.19 | 0.517 | 0.19 |

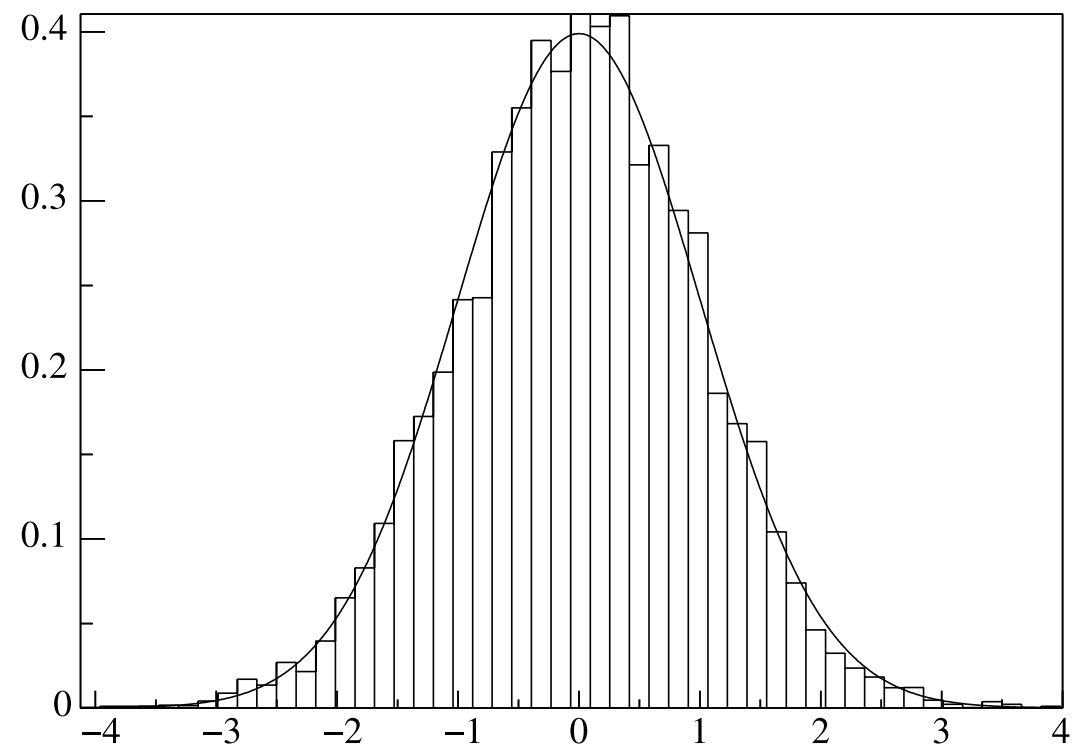

FIG. 5. *Limiting distribution of the first component of $\vartheta_n^{f,A}$ (RRIS) for the option of Table 4 with $K = 50$.*

3.2.2. *Dupire's framework.* We consider an $I$-dimensional local volatility model in which the dynamics under the risk-neutral measure of each asset $S^i$ is supposed to be given by

$$dS_t^i = S_t^i(r\,dt + \sigma(t, S_t^i)\,dW_t^i), \qquad S_0 = (S_0^1, \ldots, S_0^d),$$

where $W = (W^1, \ldots, W^I)$ is defined and generated as in the Black–Scholes framework. The local volatility function $\sigma$ we have chosen is of the form

$$\sigma(t,x) = 0.6(1.2 - e^{-0.1t}e^{-0.001(xe^{rt}-s)^2})e^{-0.05\sqrt{t}}, \tag{3.1}$$

with $s > 0$. We know that there exists a duality between the variables $(t,x)$ and $(T,K)$ in Dupire's framework. Hence, for the formula (3.1) to make sense, one should choose $s$ equal to the spot price of the underlying asset so that the bottom of the smile is located at the forward money. We refer to Figure 7 for an overview of the smile.

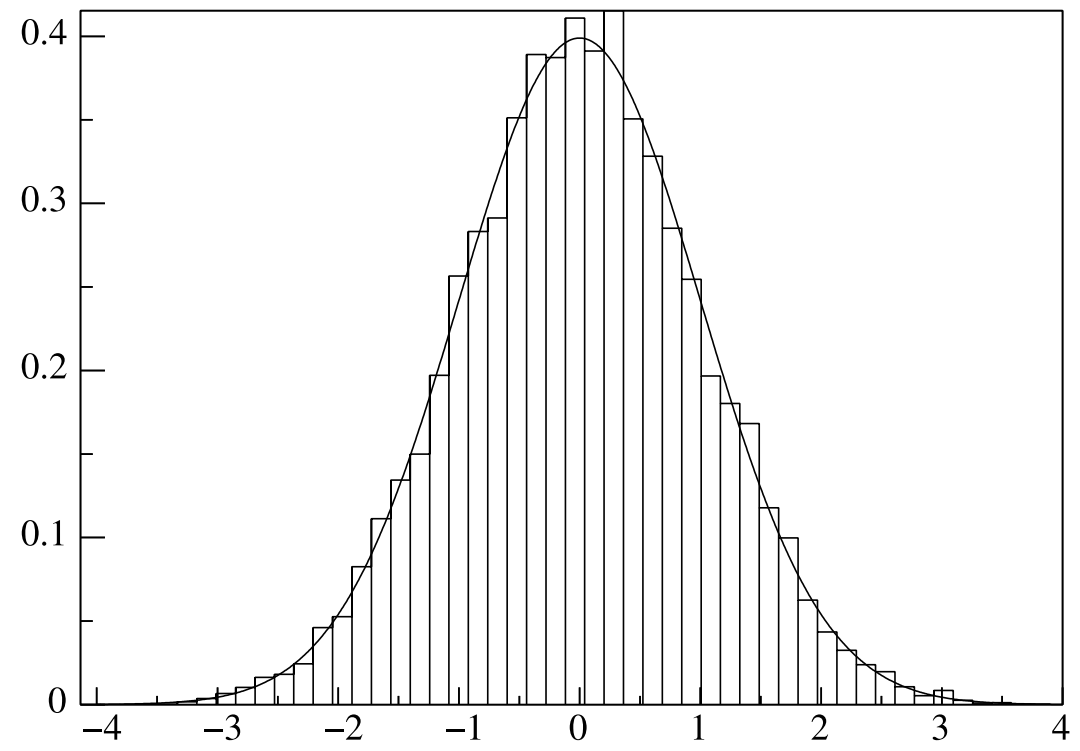

FIG. 6. *Limiting distribution of $M_n(\theta_n^{f,A}, f)$ (RRIS) for the option of Table 4 with $K = 50$.*

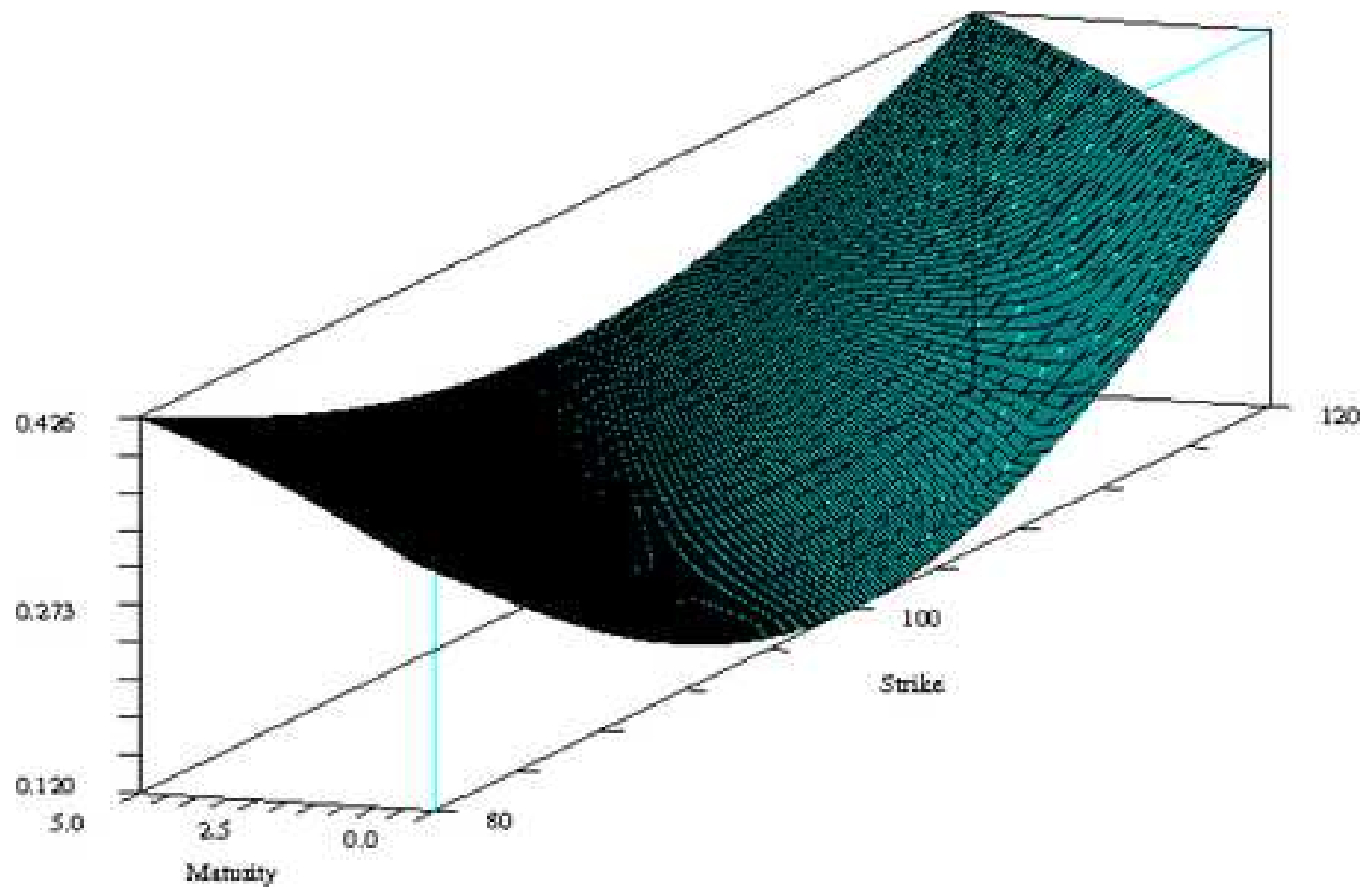

FIG. 7. *Local volatility function.*

*Best-of option.* We consider options with payoffs $(\max_{1\leq i\leq I}\omega_i S_T^i - K)_+$, where $K > 0$ and $(\omega^1, \dots, \omega^I)$ is a vector of positive weights. The payoffs belong to $\mathcal{H}_1$. To discretize the dynamics, we use an Euler scheme with $N = 100$ time steps per year. The results of Table 5 are encouraging. The RRIS algorithm with $A$ defined as in the barrier basket option case reduces the variance by 6, whereas it only increases the computational time by 2, which leads to a gain of 6/2. We do not present any results for the RIS algorithm because the extra computational time it requires makes it noncompetitive.

**Conclusion.** We propose a fully automatic adaptive importance sampling technique for the computation of $\mathbb{E}(f(G))$, where $f:\mathbb{R}^d \to \mathbb{R}$ and $G$ is a standard $d$-dimensional normal random vector. For a large class of functions $f$, including many financial payoffs, we prove that our estimator is convergent and asymptotically normal with optimal limiting variance. Note that all of the convergence results stated in Theorems 2.2, 2.3, Corollary 2.4, Propositions 2.7, 2.14, Lemma 2.8 and Remarks 2.10, 2.11, 2.16 still hold if $M_n(\theta, g)$ is defined as $\frac{1}{n}\sum_{i=1}^n g(\tilde{G}_i + \theta)e^{-\theta\cdot\tilde{G}_i - |\theta|^2/2}$ for any sequence $(\tilde{G}_i)_{i\geq 1}$ of i.i.d. $d$-dimensional standard normal random vectors and, in particular, when this sequence is independent from the one, $(G_i)_{i\geq 1}$, used to compute $(\vartheta_n^{f,A})_{n\geq 1}$. Our numerical experiments confirm the effectiveness of our estimator: in comparison with the crude Monte Carlo method, the computation time needed to achieve a given precision is divided by a factor between 3 and 15. Moreover, they suggest that the convergence and asymptotic normality of the estimator still hold under weaker assumptions on the function $f$. In view of these numerical results and the definition of $\mathcal{V}_1$, it would be natural to investigate the class of functions $f$ such that, for some constants $\lambda > 0$ and $\beta \in [0, 2)$,

$$\forall \varphi : \mathbb{R}^d \to \mathbb{R}^d \ C^\infty \text{ and vanishing outside } B(0, M),$$
$$\left| \int_{\mathbb{R}^d} f \nabla \cdot \varphi(x)\, dx \right| \leq \lambda e^{M^\beta} \|\varphi\|_\infty.$$

TABLE 5
*Best-of option in dimension* 12 *with* $\rho = 0.5$, $r = 0.05$, $T = 1$, $n = 50{,}000$ *and* $\omega^i = 1$, $S_0^i = 50$ *for all* $i = 1, \dots, I$

| $K$ | Price | Price MC | Var MC | Price RRIS | Var RRIS |
|---|---|---|---|---|---|
| 70 | 3.260 | 3.236 | 137 | 3.299 | 24.50 |
| 80 | 1.901 | 1.917 | 94.23 | 1.905 | 14.09 |
| 90 | 1.220 | 1.253 | 67.70 | 1.227 | 9.41 |

Unfortunately, we have thus far not been able to derive the asymptotic properties of our estimator for such functions. In this work, we have focused on importance sampling. A natural extension would be to investigate the coupling with stratification techniques in the spirit of [8]. In particular, it would be interesting to combine the present importance sampling algorithm with the adaptive stratified sampling methods recently proposed in [6] (adaptive optimization of the proportions of random drawings made in the different strata) and [5] (adaptive optimization of the stratification direction $e \in \mathbb{R}^d$ for a standard normal random vector when the strata are given by $\{x \in \mathbb{R}^d : e \cdot x \in [y_{i-1}, y_i)\}$ with $-\infty = y_0 < y_1 < y_2 < \cdots < y_I = +\infty$).

Université Paris-Est
6 et 8 avenue Blaise Pascal
77455 Marne La Vallée
Cedex 2
France
E-mail: jourdain@cermics.enpc.fr

Ecole Nationale Supérieure
de Techniques Avancées
42 bd Victor 75015
Paris
France
E-mail: jerome.lelong@ensta.fr